\documentclass[10pt]{amsart}
\numberwithin{equation}{section}
\usepackage{pdfsync}
\usepackage{graphicx}
\usepackage{enumerate}
\usepackage{xspace}
\usepackage{boxedminipage}
\usepackage{subfigure}
\usepackage{hyperref} 
\usepackage{epsfig}
\usepackage{latexsym}
\usepackage{amsmath,amssymb,amsxtra,amsthm}
\usepackage{stmaryrd}
\usepackage{geometry}                % See geometry.jpg to learn the layout options. There are lots.
\usepackage{commath}
\usepackage{mathrsfs} %nice calligraphic characters
\usepackage{ifthen}
\provideboolean{showchanges}
\setboolean{showchanges}{false}%%false
\newcommand{\changes}[1]{
  \ifthenelse{\boolean{showchanges}} {{\bl{#1}}} {#1}
}
\usepackage{multicol}        % used for the two-column index
\usepackage[bottom]{footmisc}% places footnotes at page bottom

\makeindex             % used for the subject index
\newcommand{\Proofname}{Proof}
 
\newtheoremstyle%
    {plain}% name
    {}% Space before, vuoto = `valore di default'
    {}% Space after
    {\mdseries\slshape}% body font
    {}% Indent (empty = no indent, \parindent = para indent)
    {\bfseries}% Thm head font
    {.}% Punctuation after the heading
    {1.0ex}% Space after heading: \newline = to start at next line
    {}% Thm head spec (can be left empty, meaning `normal')
    
\newtheoremstyle
    {note}% name
    {}% Space before, vuoto = `valore di default'
    {}% Space after
    {}% body font
    {}% Indent (empty = no indent, \parindent = para indent)
    {\bfseries}% Thm head font
    {.}% Punctuation after the heading
    {1.0ex}% Space after heading: \newline = to start at next line
    {}% Thm head spec (can be left empty, meaning `normal')
    
\swapnumbers
\newtheorem{The}[subsection]{Theorem}

\newtheorem{Prop}[subsection]{Proposition}

\theoremstyle{note}

\newtheorem{Assum}[subsection]{Assumption}

\renewcommand{\vec}[1]{\ensuremath{\boldsymbol{#1}}}
\newcommand{\myall}{\ensuremath{\quad \forall}}
\newcommand{\pdt}{\ensuremath{\partial_t}}

\newcommand{\Reals}{\ensuremath{{\mathbb{R}}}}

\renewcommand{\O}{\ensuremath{{\Omega}}}
\newcommand{\Ot}{\ensuremath{{\Omega_t}}}

\newcommand{\Oc}{\ensuremath{{\hat{\Omega}}}}
\newcommand{\normal}{\ensuremath{{\vec{\nu}}}}

\newcommand{\transpose}{\ensuremath{{\mathsf{T}}}}
\newcommand{\Lp}[1]{\ensuremath{\operatorname{L}_{#1}}}

\newcommand{\Hil}[1]{\ensuremath{\operatorname{H}^{#1}}}
\newcommand{\Vc}{\ensuremath{\hat{\mathbb{V}}}}
\newcommand{\ltwop}[3]{\ensuremath{\left\langle#1,#2\right\rangle}_{#3}}
\newcommand{\dualp}[3]{\ensuremath{\left\langle#1\,\vert \, #2\right\rangle}_{#3}}
\newcommand{\ltwon}[2]{\ensuremath{\left\|#1\right\|}_{\Lp{2}#2}}

\newcommand{\linfn}[2]{\ensuremath{\left\|#1\right\|}_{\Lp{\infty}#2}}

\newcommand{\A}{\ensuremath{\vec{\mathcal{A}}}}
\newcommand{\lv}{\ensuremath{\left\vert}}
\newcommand{\rv}{\ensuremath{\right\vert}}
\usepackage{upgreek}
\newcommand{\lap}{\ensuremath{\Updelta}}

\newcommand{\T}{\ensuremath{{\mathscr{T}}}}
\newcommand{\Tc}{\ensuremath{\hat{\mathscr{T}}}}
\newcommand{\s}{\ensuremath{{{s}}}}

\newcommand{\E}{\ensuremath{{{E}}}}
\newcommand{\interior}{\ensuremath{{{i}}}}

\newcommand{\fim}{for $i=1,\dotsc,m$, }
\newcommand{\jump}[1]{\ensuremath{\left\llbracket{#1}\right\rrbracket}}
\newcommand{\estimator}{\ensuremath{\hat{\eta}}}
\usepackage{color}
\newcommand{ \bl}{\color{blue}}

\definecolor{MyGreen}{rgb} {0.05,0.4,0.05}
\definecolor{RedViolet}{rgb} {0.1,0.1,0.75}

\usepackage{algorithm, algorithmic}

\begin{document}
\bibliographystyle{unsrt}
\title[AFEM for RDSs on growing domains]{Adaptive finite elements for semilinear reaction-diffusion systems on growing domains}
%for an abbreviated version of
\author{Chandrasekhar Venkataraman}
\author{Omar Lakkis} 
\author{Anotida Madzvamuse}
\address[Chandrasekhar Venkataraman\\Omar Lakkis\\Anotida Madzvamuse]{Department of Mathematics, University of Sussex, Pevensey III, Brighton, BN1 9QH.} 
 \email[Chandrasekhar Venkataraman]{C.Venkataraman@sussex.ac.uk}\email[Omar Lakkis]{O.Lakkis@sussex.ac.uk}\email[Anotida Madzvamuse]{A.Madzvamuse@sussex.ac.uk}
\thanks{The research of CV has been supported by the EPSRC, Grant EP/G010404.}
\thanks{A version of this article has been published in Numerical Mathematics and Advanced Applications 2011 \cite{venkataraman2013adaptive}}
\maketitle
%\vspace{-10em}
\begin{abstract}
We propose an adaptive finite element method to approximate the solutions to reaction-diffusion systems on time-dependent domains and surfaces. We derive a computable error estimator that provides an upper bound for the error in the semidiscrete (space) scheme.  We reconcile our theoretical results with benchmark computations.
\end{abstract}
\maketitle
%\vspace{-2em}
\section{Introduction}\label{C:fem_apriori:secn.:setup}
%\vspace{-1em}
%Model Problem
%%%%%%%%%%%%%%%%%%%%%%%%%%%%%%%%%%%%%%%%%%%%%%%%%%%%%%%%%%%%%%%%%%%%%%%%%%%%%%%%%%%%%%%%%%%%%%%%%%%%%%%%%%%%%%%%%%%%%%%%%%%%%%%%%%%%%%%%%%%%%%%%%%%%%%%%%%%%%%%%%%%%%%%%%%%%%%%%%%%%%%%%%%%%%%%%%%%%%%%%%%%%%%%%%%%%%%%%%%%%%%%%%%%%%%%%%%%%%%%%%%%%%%%%%%%%%%%%%%%%%%%%%%%%%%%%%%%%%%%%%%%%%%%%%%%%%%%%%
Our model problem consists of a system of chemicals that  are coupled only through the reaction terms and diffuse independently of each other.  Given an integer $m\geq1$, let $\vec{u}\left(\vec{x},t\right)$ be an ($m\times1$) vector of
concentrations of chemical species, with $\vec{x} \in\Ot\subset
\mathbb{R}^2$,  the spatial variable and  $t \in [0,T],
\ T>0,$ the time variable. The model we shall consider is of the following form (see \cite{ano2000} for details of the
derivation): 
 find $u_i$, functions from $\Ot$ into $\Reals$, such that \fim $u_i$ satisfies  \begin{equation}\label{C:derivation:eqn:model_problem}
 \begin{split}
\begin{cases}
 \pdt{u}_i(\vec{x},t)-D_i\lap{u}_i(\vec{x},t)+\nabla\cdot\left[\vec{a}u_i\right](\vec{x},t)=f_{i}\left(\vec{u}(\vec{x},t)\right), &\vec{x}\in\Ot,t\in(0,T],\\
[\normal\cdot\nabla{u}_i](\vec{x},t)=0, &\text{ } \vec{x}\in\partial\Ot, t>0,\\
u_i(\vec{x},0)=u_i^{0}(\vec{x}), &\text{ }\vec{x}\in\O_0,
\end{cases}
\end{split}
\end{equation}
where $\Omega_{t}$ is a simply connected bounded continuously deforming domain with respect to
$t$, with  Lipschitz boundary $\partial\Omega_{t}$ at time 
$t \in [0,T]$. The vector of nonlinear coupling terms $\vec{f}:=({f}_1,\dotsc,{f}_m)^\transpose$ is assumed to be locally Lipschitz-continuous,
$\vec{D}:=({D}_1,\dotsc,{D}_m)^\transpose$ is a vector of strictly positive diffusion coefficients,  $\vec{a}=(a_1,\dotsc,a_d)^\transpose$ is a flow velocity generated by the evolution of the domain and the initial data $\vec{u}^{0}(\vec{x})$ is a bounded vector valued function.  Systems of this form arise  in the theory of biological pattern formation \cite{venkataramanthesis}.

Let $\Oc$ be a simply connected  time-independent reference domain with Lipschitz boundary.   We assume there exists a time-differentiable family of $C^\infty$-diffeomorphisms  $\vec{\mathcal{A}} :\Oc \times[0,T] \rightarrow\Omega_t$ such that at each instant $t\in[0,T]$ and for each $\vec x\in\Omega_t$ there exists a $\vec\xi\in\hat\Omega$ such that
\begin{equation}
\label{C:derivation:eqn:mapping_def}
\A(\vec{\xi},t)=\vec{x}. 
\end{equation}

Based on the derivation presented in \cite{venk10fem} 
%we write  an equivalent problem posed on the {\it time independent} reference domain:
%For $i=1,\dotsc,m$, $\hat{u}_i(\vec{\xi},t)$,  the pullback of the function ${u}_i(\A(\vec \xi, t),t)$, i.e.
%$\hat{u}(\vec \xi, t)=u(\A(\vec \xi, t),t)$, satisfies the following semilinear reaction-diffusion-convection system posed on the reference domain: 
%\begin{equation}
%\begin{split}
%\label{C:derivation:eqn:reference_model_problem}
%\pdt{\hat{u}}_i\left(\vec{\xi},t\right)+\hat{u}_i\left(\vec{\xi},t\right)\nabla\cdot\vec{a}\left(\A(\vec{\xi},t),t\right)
%={f}_i\left(\hat{\vec{u}}\left(\vec{\xi},t\right)\right)&+D_{i}\left[\Dif\A^{-1}:\Dif\A^{-1}\right]\left(\A\left(\vec{\xi},t\right)\right)\lap\hat{u}_i\left(\vec{\xi},t\right)\\
%+D_{i}\bigg(\big[\lap\mathcal{A}^{-1}_1,\dotsc,\lap\mathcal{A}^{-1}_d\big]&\left(\A\left(\vec{\xi},t\right)\right)\bigg)\nabla\hat{u}_i\left(\vec{\xi},t\right)\quad\vec{\xi} \in \Oc, t \in (0,T],\\
%[\hat \normal\cdot\vec K\nabla{\hat u}_i](\vec{\xi},t)=&0, \quad \vec{\xi}\in\partial\Oc, t>0,\\
%\hat u_i(\vec{\xi},0)=&\hat u_i^{0}(\vec{x}), \quad \vec{\xi}\in\Oc,
%\end{split}
%\end{equation}
% where $\vec{K}$ is the  Jacobian of the inverse diffeomorphism $\A^{-1}$.  
%
%
%
%For the purposes of constructing a finite element discretisation,
 we introduce a {\em weak formulation} associated with Problem (\ref{C:derivation:eqn:model_problem}) on the {\it  time-independent reference domain}. The problem is to find
$\hat{u}_i \in \Lp{2}{\big(0,T;{\Hil{1}{(\smash{\Oc})}}\big)}$ with $\pdt\hat{u}_i\in\Lp{2}{\big(0,T;\changes{\Hil{1}{(\smash{\Oc})}^{\prime}}\big)}$ such that  
\changes{for all $t\in(0,T]$},
\begin{equation}
\begin{split}
\label{C:fem_apriori:eqn:cont_weak_form_reference_ddtju}
\ltwop{\pdt(J\hat{u}_i)}{\hat{\chi}}{\Oc}+\ltwop{{D}_iJ\vec{K}{\nabla}\hat{u}_i}{\vec{K}\nabla\hat{\chi}}{\Oc}&=\ltwop{Jf_i(\hat{\vec{u}})}{\hat{\chi}}{\Oc},  \myall \hat{\chi} \in \Hil{1}{(\smash{\Oc})}.
\end{split}
\end{equation}
Here  \changes{$\Hil{1}{(\smash{\Oc})}^{\prime}$ is the dual of $\Hil{1}{(\smash{\Oc})}$ equipped with the norm
\begin{equation}
\| v\|_{\Hil{1}{(\Oc)}^{\prime}}:=\sup_{w\in\Hil{1}{(\Oc)},w\neq0}\frac{\dualp{v}{w}{\Hil{1}{(\Oc)}^{\prime}\times\Hil{1}{(\Oc)}}}{\|w\|_{\Hil{1}{(\Oc)}}}.
\end{equation}
The matrix $\vec K$ and $J$ are the inverse and determinant of the Jacobian of the diffeomorphism $\A$ respectively.}  
%\vspace{-2em}

\section{A posteriori error estimates}\label{C:fem_apriori:secn.:aposteriori}
%\vspace{-1em}
Here we state a Theorem, and the associated Assumptions under which the Theorem holds, that shows the error in the {\it semidiscrete} scheme  can be bounded by a computable  a posteriori error estimator,   based on the element {\em residual}.
Our strategy to derive an a posteriori error estimate is similar to that employed in \cite{kruger2003posteriori}. We use energy arguments to show the residual is an upper bound for the error and the analysis is similar to the a priori case we have considered elsewhere \cite{venk10fem}. For the details of the proofs we refer to \cite{venkataramanthesis}.

 We start by stating the semidiscrete scheme, find $\hat{u}^h_i:[0,T]\rightarrow\Vc$, such that for $i=1,\dotsc,m$,
%%%%%%%%%%%%%%%%%%%%%%%%%%%%%%%%%%%%%%%%%%%%%%%%%%%%%%%%%%%%%%%%%%%%%%%%%%%%%
\begin{equation}
\begin{cases}
\begin{aligned}
\label{C:fem_apriori:eqn:sd_weak_form_aposteriori}
\ltwop{\pdt(J\hat{u}^h_i)}{\hat{\Phi}}{\Oc}+\ltwop{{D}_iJ\vec{K}{\nabla}\hat{u}^h_i}{\vec{K}\nabla\hat{\Phi}}{\Oc}&=\ltwop{J\widetilde{f}_i(\hat{\vec{u}}^h)}{\hat{\Phi}}{\Oc} \myall \hat{\Phi} \in \Vc\changes{\text{ and }t\in(0,T]},\\
\hat{u}^h_i(0)&=\Lambda^h\hat{u}_i^0, 
\end{aligned}
\end{cases}
\end{equation}
%%%%%%%%%%%%%%%%%%%%%%%%%%%%%%%%%%%%%%%%%%%%%%%%%%%%%%%%%%%%%%%%%%%%%%%%%%%%%
where $\Vc$  is a standard FE space made up of piecewise polynomial functions and $\Lambda^h : H^1(\smash{\Oc}) \rightarrow \Vc$ is the Lagrange interpolant.
\begin{Assum}[Applicability of the MVT]\label{assum:mvt}
We assume that 
 \begin{equation}
 \linfn{\vec{f}^{\prime}}{\left(\operatorname{dom}(\hat{\vec{u}})\right)}+\linfn{\vec{f}^{\prime}}{\left(\operatorname{dom}(\hat{\vec{u}}^h)\right)}<\tilde{C}.
\end{equation}
Note this assumption is satisfied if we assume a global smallness condition on the mesh-size \cite{venk10fem} and that the continuous problem is well posed \cite{venkataraman2010global}.
\end{Assum}
%%%%%%%%%%%%%%%%%%%%%%%%%%%%%%%%%%%%%%%%%%%%%%%%%%%%%%%%%%%%%%%%%%%%%%%%%%%%%
We define the error in the semidiscrete scheme  
\begin{equation}\label{C:fem_apriori:eqn:error_defn}
\hat{\vec{e}}(t):=\hat{\vec{u}}^h(t)-\hat{\vec{u}}(t), \text{ for } t\in[0,T].
\end{equation}

\begin{Assum}[Dominant energy norm error]\label{assum:energy_l2_dominance}
Since we are primarily interested in problems posed on long time intervals, we wish to circumvent the use of  Gronwall's inequality.  To this end we assume that the error in the $\Lp{2}(0,T;\Lp{2}(\Oc))$ norm converges faster than the  error in the $\Lp{2}(0,T;\Hil{1}(\Oc))$ norm. We  assume there exists $C^\dagger,C> 0$ and $r\in(0,1]$ independent of the mesh-size $\hat{h}$ such that
\begin{equation}\label{eqn:dominant-energy-evolving}
\int_0^T\ltwon{\hat{\vec{e}}}{(\Oc)^m}^2\leq{C^\dagger}\hat{h}^{2r}\sum_{i=1}^m\int_0^T\ltwon{\nabla\hat{e}_i}{(\Oc)}^2,
\end{equation}
thus
\begin{equation}
\int_0^T\ltwon{{\vec{e}}}{(\Ot)^m}^2\leq{C}\hat{h}^{2r}\sum_{i=1}^m\int_0^T\ltwon{\nabla{e}_i}{(\Ot)}^2,
\end{equation}
where we have used the equivalence of norms between the reference and evolving domains. 
\end{Assum}

We note assumptions of this type have been used previously in \cite{kruger2003posteriori} and  \cite{medina1996error} to obtain a posteriori estimates for quasilinear reaction-diffusion  and nonlinear convection-diffusion problems.

We start by introducing the {\em residual} $ \hat{R}_i\in\Hil{1}(\Oc)^{\prime}$ (the dual of $\Hil{1}{(\Oc)}$) a.e. in $[0,T]$ which satisfies 
\begin{equation}\label{C:apriori_fem:eqn:residual_definition}
\dualp{\hat{R}_i}{\hat{\chi}}{}:=\dualp{\pdt({J\hat{u}^h}_i)-D_i\nabla\cdot\left(J\vec{KK}^\transpose\nabla{\hat{u}^h_i}\right)-Jf_i(\hat{\vec{u}}^h)}{\hat{\chi}}{} \myall\hat{\chi}\in\Hil{1}(\Oc),
\end{equation}
where $\dualp{\cdot}{\cdot}{}$ denotes the duality pairing between $\Hil{1}{}$ and its dual. We now show the residual is an upper bound for the error.
\begin{Prop}[Upper bound for the error]
Suppose Assumptions \ref{assum:mvt} and
\ref{assum:energy_l2_dominance} hold. Let $\hat{\vec{u}}$ satisfy (\ref{C:fem_apriori:eqn:cont_weak_form_reference_ddtju}) and let the error $\hat{e}_i$ and the residual $\hat{R}_i$ be as in (\ref{C:fem_apriori:eqn:error_defn}) and (\ref{C:apriori_fem:eqn:residual_definition}) respectively.
 If the mesh-size satisfies  a smallness condition   (see \cite{venkataramanthesis} for details), then
 \begin{equation}\label{C:fem_apriori:eqn:error_residual_upper_bound}
 \begin{split} 
\ltwon{{\vec{e}}(T)}{(\Ot)^m}^2+\sum_{i=1}^mD_i\int_0^T\ltwon{\nabla{e}_i}{(\Ot)}^2\leq\ltwon{{\vec{e}}(0)}{(\O_0)^m}^2&+2\sum_{i=1}^m\int_0^T\dualp{\hat{R}_i}{\hat{e}_i}{}.
\end{split}
\end{equation}
\end{Prop}
We now introduce a concrete error estimator. For simplicity we restrict the discussion to the case of $\mathbb{P}^1$ elements and regular triangulations, the results may be straightforwardly generalised to higher order elements.
For any simplex $\s $ of the triangulation $\hat{\T}$ we denote by $\hat{h}_\s$ the \changes{diameter} of $\s$. Let $\E_{\s }$ be the set of three edges of $\s $. Let $\E_{\interior}$ be an edge on the interior of $\Oc$, with outward pointing (with respect to $\s$) normal  $\normal$. We denote by $\jump{\nabla\phi\cdot\normal}_{\E_{\interior}}$ the jump of $\nabla\phi\cdot\normal$ across the edge $\E_{\interior}$. For boundary edges we take $\jump{\nabla\phi\cdot\normal}=2\nabla\phi\cdot\normal$. The local error indicator  is given by
\begin{equation}\label{C:fem_apriori:eqn:estimator_defn}
\begin{split}
( \estimator_{i\vert{\s}})^2
:=\hat{h}_{\s}^2&\ltwon{\pdt(J\hat{u}^h_i)-D_i\nabla\cdot\left(J\vec{KK}^\transpose\nabla{\hat{u}^h_i}\right)-Jf_i(\hat{\vec{u}}^h)}{(\s)}^2\\
&+\frac{1}{2}\sum_{e\in\E_{\s}}\lv{e}\rv\ltwon{D_i\jump{J\vec{KK}^\transpose\nabla\hat{u}^h_i\cdot\normal}}{(e)}^2.
\end{split}
\end{equation}
\begin{Prop}[Residual bound]
Let  ${R}_i$ and $\estimator_i$, $i=1,\dotsc,m$ be defined by (\ref{C:apriori_fem:eqn:residual_definition}) and (\ref{C:fem_apriori:eqn:estimator_defn}) respectively. There exists a $C>0$ that depends only on the shape regularity of the triangulation $\hat{\T}$ such that \fim
\begin{equation}\label{C:apriori_fem:eqn:residual_bound}
\Bigg\vert\int_0^T\dualp{\hat{R}_i}{\hat{\chi}}{}\Bigg\vert
\quad\leq{C}\int_0^T\left(\sum_{\s\in\hat{\T}}(\estimator_{i\vert{\s}})^2\right)^{1/2}\int_0^T\changes{\Big\|{\hat{\chi}}\Big\|_{\Hil{1}{(\Oc)}}} \myall \hat{\chi}\in\changes{\Lp{2}([0,T];\Hil{1}{(\Oc)})}.
\end{equation}
\end{Prop}

To complete the bound  of the error by the estimator, we make an assumption about the error in the approximation of the initial data.
\begin{Assum}[Dominated initial error]\label{assum:initial_l2_dominance}
We assume that the initial error in the $\Lp{2}(\Oc)$ norm converges faster than the  error in the $\Lp{2}(0,T;\Hil{1}(\Oc))$.  We assume there exists $C> 0$ and $r\in(0,1]$ both independent of the mesh-size $\hat{h}$ such that
\begin{equation}
\ltwon{\hat{\vec{e}}(0)}{(\Oc)}^2\leq{C}\hat{h}^{2r}\sum_{i=1}^m\int_0^T\ltwon{\nabla\hat{e}_i}{(\Oc)}^2.
\end{equation}
\end{Assum}
\begin{The}[A posteriori error estimate for the semidiscrete scheme]\label{C:fem_apriori:Thm:aposteriori}
Let Assumptions \ref{assum:mvt}, \ref{assum:energy_l2_dominance} and \ref{assum:initial_l2_dominance} hold. Let the error $\hat{e}_i$ and  the estimator $\estimator_i,$ $i=1,\dotsc,m$ be defined by  (\ref{C:fem_apriori:eqn:error_defn}) and (\ref{C:fem_apriori:eqn:estimator_defn}) respectively.
 If the mesh-size is sufficiently small, for some $C>0$, we have
 \begin{equation}\label{C:fem_apriori:Thm:eqn:aposteriori}
 \sum_{i=1}^mD_i\int_0^T\ltwon{\nabla{e}_i}{(\Ot)}^2\leq\sum_{{i=1}}^{m}{C}\int_{0}^T\sum_{\s\in\hat{\T}}(\estimator_{i\vert\s})^2.
 \end{equation}
 \end{The}
 
Since the estimator $\vec{\estimator}$ is an upper bound for the error,  we use it to drive a space-adaptive scheme.
To ensure the efficiency of the adaptive scheme, we would have to show the estimator was also a lower bound for the error and we leave this extension for future work.
%%%%%%%%%%%%%%%%%%%%%%%%%%%%%%%%%%%%%%%%%%%%%%%%%%%%%%%%
%\vspace{-2em}
\section{Numerical results}
%\vspace{-1em}
Here we reconcile our theoretical results with numerical computations. We start by presenting a time discretisation of the {\it semidiscrete} scheme \changes{(\ref{C:fem_apriori:eqn:sd_weak_form_aposteriori})}, we discretise in time using a modified implicit Euler method \cite{madzvamuse2007modified}, in which the reaction terms are treated semi-implicitly while the diffusive terms are treated fully implicitly: find $(\hat{U}^h_i)^n\in\Vc^n$, such that for $i=1,\dotsc,m$,  $n=1,\dotsc,N$,
%%%%%%%%%%%%%%%%%%%%%%%%%%%%%%%%%%%%%%%%%%%%%%%%%%%%%%%%%%%%%%%%%%%%%%%%%%%%%
\begin{equation}
\begin{cases}
\begin{aligned}
\label{C:fem_apriori:eqn:fd_weak_form_aposteriori}
\ltwop{\frac{1}{\tau}\bar{\partial}\left[J\hat({U}^h_i)\right]^n}{\hat{\Phi}^n}{\Oc}+&\ltwop{{D}_i[J\vec{K}{\nabla}\hat({U}^h_i)]^n}{[\vec{K}\nabla\hat{\Phi]^n}}{\Oc}\\
&=\ltwop{J^n\widetilde{f}_i((\hat{U}^h_i)^n,(\hat{U}^h_i)^{n-1})}{\hat{\Phi}^n}{\Oc} \myall \hat{\Phi} \in \Vc^n,\\
(\hat{U}^h_i)^0&=\Lambda^h\hat{u}_i^0.
\end{aligned}
\end{cases}
\end{equation}

The adaptive algorithm we consider is based on the
equidistribution marking strategy \cite[Alg.~1.19, pg. 45]{schmidt2005design}, where elements are marked for refinement and coarsening with the goal of equidistributing the estimator value over all mesh elements.  The marking strategy takes two parameters: the tolerance of the adaptive algorithm $tol$ and a parameter  $\theta\in(0,1)$. At each timestep elements are marked for refinement according to the following algorithm:
\newcommand{\Algoname}[1]{\ensuremath{\text{\textsf{#1}\xspace}}}
\subsection{$\Algoname{Equidistribution strategy (refinement)}$} 
\begin{algorithmic}
\STATE Start with $\Tc^n_0$ the initial triangulation at time $n$, tolerance $tol$ and parameter $\theta$
\STATE $k:=0$
\STATE solve the discrete linear problem on the mesh $\Tc^n_k$
\STATE compute global error estimator $\vec{{\estimator}}$ and local error indicators $\vec{\estimator}_{\vert\s}$ 
\WHILE{$\vec{\estimator}>tol$}
\FORALL{$\s\in\Tc^n_k$}
\IF{$\vec{\estimator}_{\vert\s}>\theta*tol/N$ \COMMENT{where $N$ is the number of \changes{elements} of the triangulation}
}
\STATE mark $\s$ for refinement \COMMENT{elements are also marked for coarsening at this stage}
\ENDIF
\ENDFOR
\STATE adapt mesh $\Tc^n_k$ to give $\Tc^n_{k+1}$
\STATE $k:=k+1$
\STATE solve the discrete linear problem on the mesh $\Tc^n_k$
\STATE Compute global error estimator $\vec{{\estimator}}$ and local error indicators $\vec{\estimator}_{\vert\s}$
\ENDWHILE 
\end{algorithmic}

Elements are marked for coarsening in a  similar way to the above, the difference being that if the local error indicator plus a coarsening indicator is less than a given tolerance on an element then the element is marked for coarsening \cite[pg. 48]{schmidt2005design}.

\begin{figure}[]
\begin{minipage}{0.42\linewidth}
\centering
\includegraphics[trim = 0mm 0mm 0mm 0mm,  clip,  width=\linewidth]{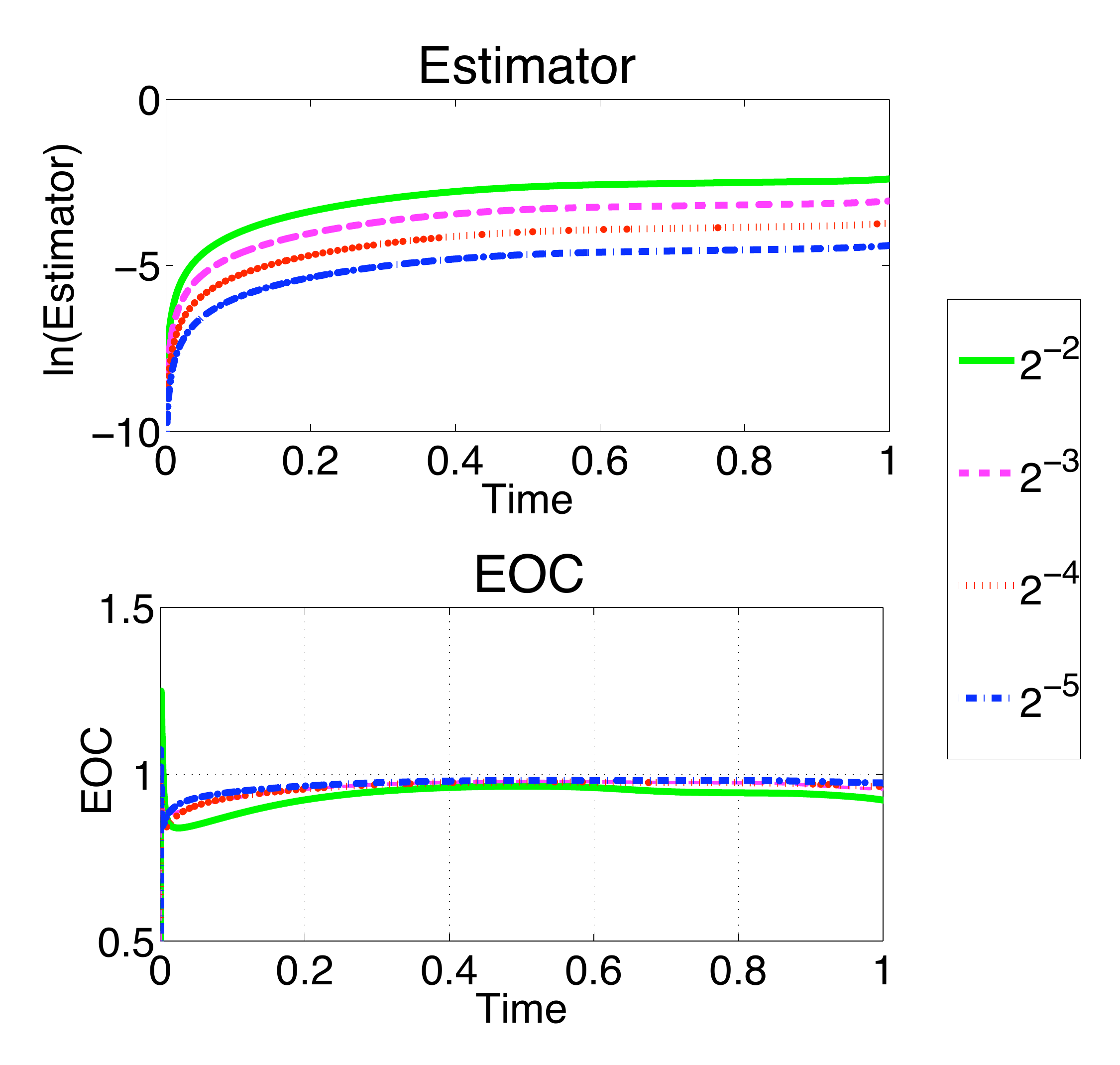} %
\caption[Values, EOC and effectivity index of the residual estimator ]{The log of the $\Lp{2}{\left([0,T]\right)}$ norm of the estimator $\estimator$ (cf. (\ref{C:fem_apriori:eqn:estimator_defn})) and  the EOC of the estimator  against time for a benchmark computation. The legend  indicates the mesh-size $\hat{h}$  for each simulation.}\label{C:fem_apriori:fig:residual_estimator}
\end{minipage}
\hspace{0.5em}
%%%%%%%%%%%%%%%%%%%%%%%%%%%%%%%%%%%%%%%%%%%%%%%%%%%%%%%
%DOF UHCH space adaptive
%%%%%%%%%%%%%%%%%%%%%%%%%%%%%%%%%%%%%%%%%%%%%%%%%%%%%%%%
\begin{minipage}{0.52\linewidth}
\centering
\includegraphics[trim = 0mm 0mm 0mm 7mm,  clip,  width=0.8\linewidth]{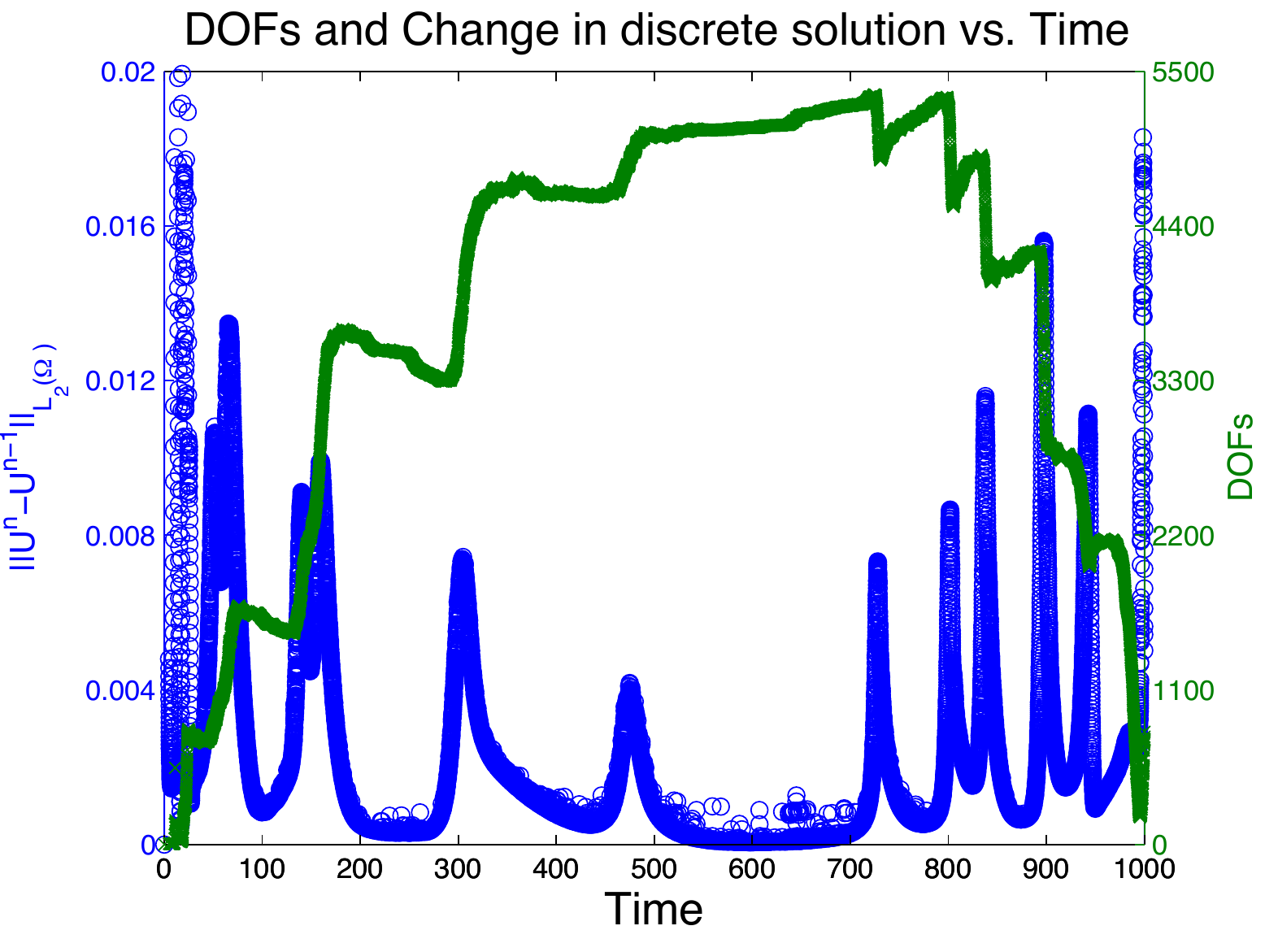}
\caption{The number of DOFs  (green crosses) and the change in discrete solution (blue circles) vs. time for the Schnakenberg kinetics on a domain with evolution of the form (\ref{eqn:SBLI_evolution}). The number of DOFs appears positively correlated with the domain size.  Bifurcations in the discrete solution correspond to spikes in the change in discrete solution. Spot-splitting bifurcations lead to increases in DOFs, while spot-annihilation or -merging results in decreases in DOFs.}
\label{C:numerics:fig:periodic_DOF}
\end{minipage}
\end{figure}

For numerical testing, we consider Problem (\ref{C:derivation:eqn:model_problem}) equipped with the Schnakenberg kinetics
\cite{prigo68}:
 \begin{equation}\label{eqn:Schnakenberg_RDS}
 f_1(\vec u)=\gamma\left(k_1-u_1+u_1^2u_2\right) \mbox{ and }
 f_2(\vec u)=\gamma\left(k_2-u_1^2u_2\right),
  \end{equation} 
  \changes{where $0<\gamma,k_1,k_2<\infty$}.
The details of the implementation of the scheme are described elsewhere \cite{venk10fem}. We consider a simulation of an RDS equipped with the Schnakenberg kinetics, with parameter values $\vec D=(1,10)^\transpose$,  $k_1=0.1$,  $k_2=0.9$,  $\gamma=1$, adding source term such that the exact solution is known, on a domain with evolution of the form
\begin{equation}
\A(\vec \xi,t)=\vec \xi(1+\sin(\pi t)), \quad\vec \xi\in[0,1]^2,t\in[0,1]. 
\end{equation}
We used a sufficiently small timestep such that the error due to the time discretisation is negligible.
The estimator values and EOC for a series of refinements is plotted in Figure \ref{C:fem_apriori:fig:residual_estimator}  and we observe an EOC of 1,  as expected for $\mathbb{P}^1$ elements providing numerical evidence for Theorem \ref{C:fem_apriori:Thm:aposteriori}.  

\begin{figure}[]
\centering
\subfigure[][$t=50$]{
\includegraphics[trim = 30mm 185mm 20mm 30mm,  clip, width=0.3\textwidth]{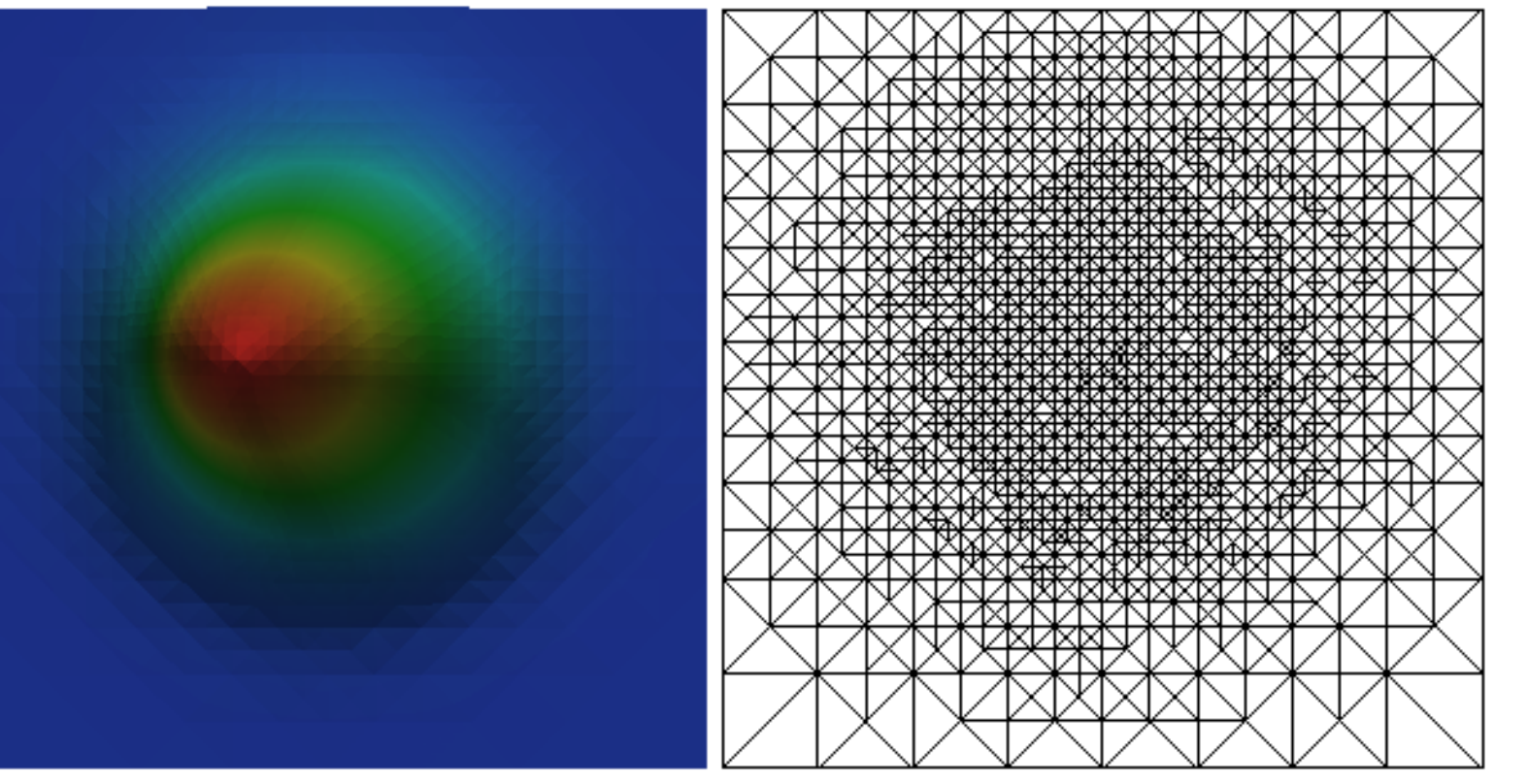}
\label{C:numerics:fig:uniform_t_50}
}
\centering
\subfigure[ ][$t=160$]{
\includegraphics[trim = 25mm 185mm 20mm 30mm,  clip,  width=0.3\textwidth]{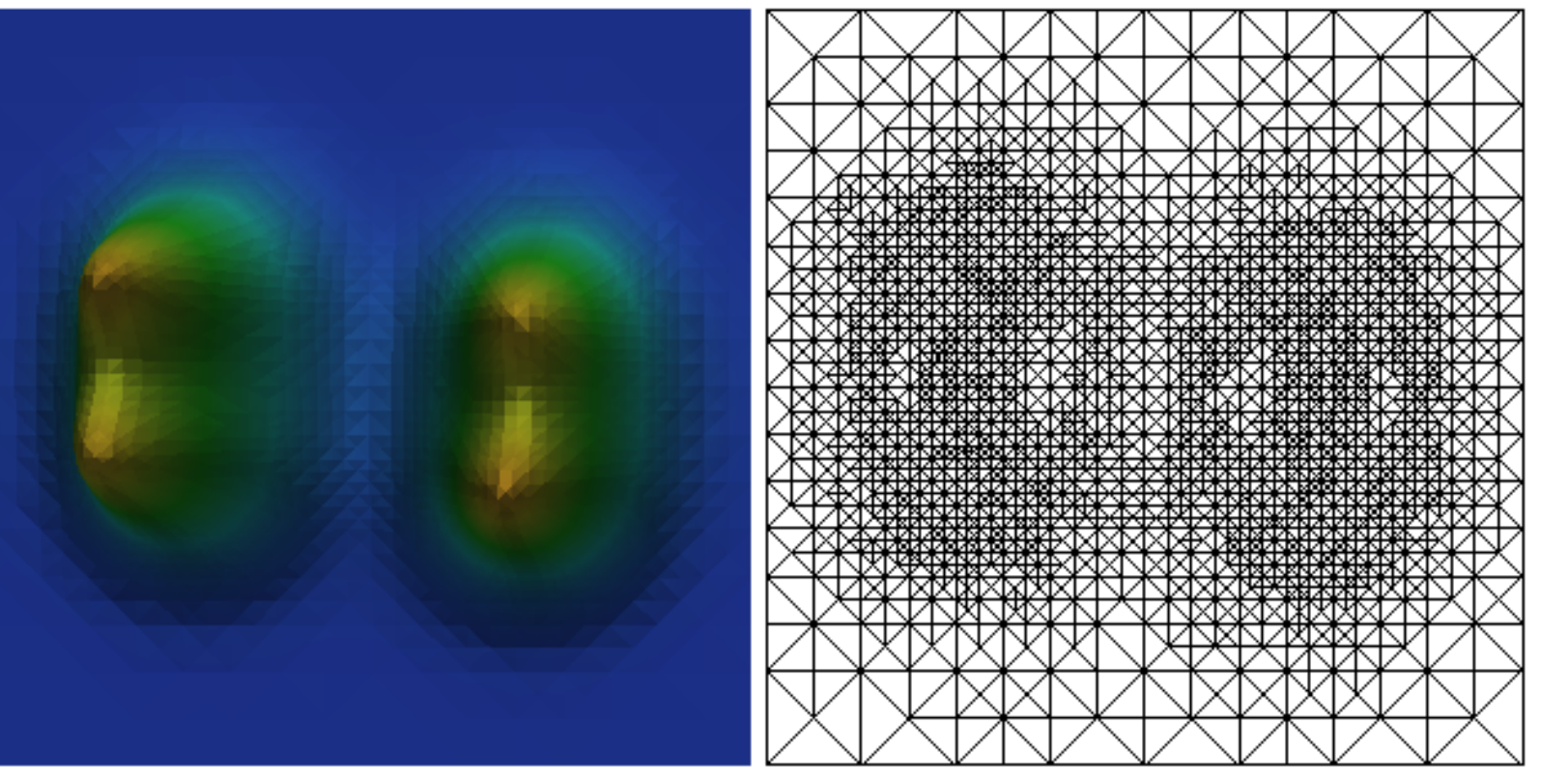}
\label{C:numerics:fig:uniform_t_160}}
\centering
\subfigure[ ][$t=380$]{
\includegraphics[trim = 25mm 185mm 20mm 30mm,  clip, width=0.3\textwidth]{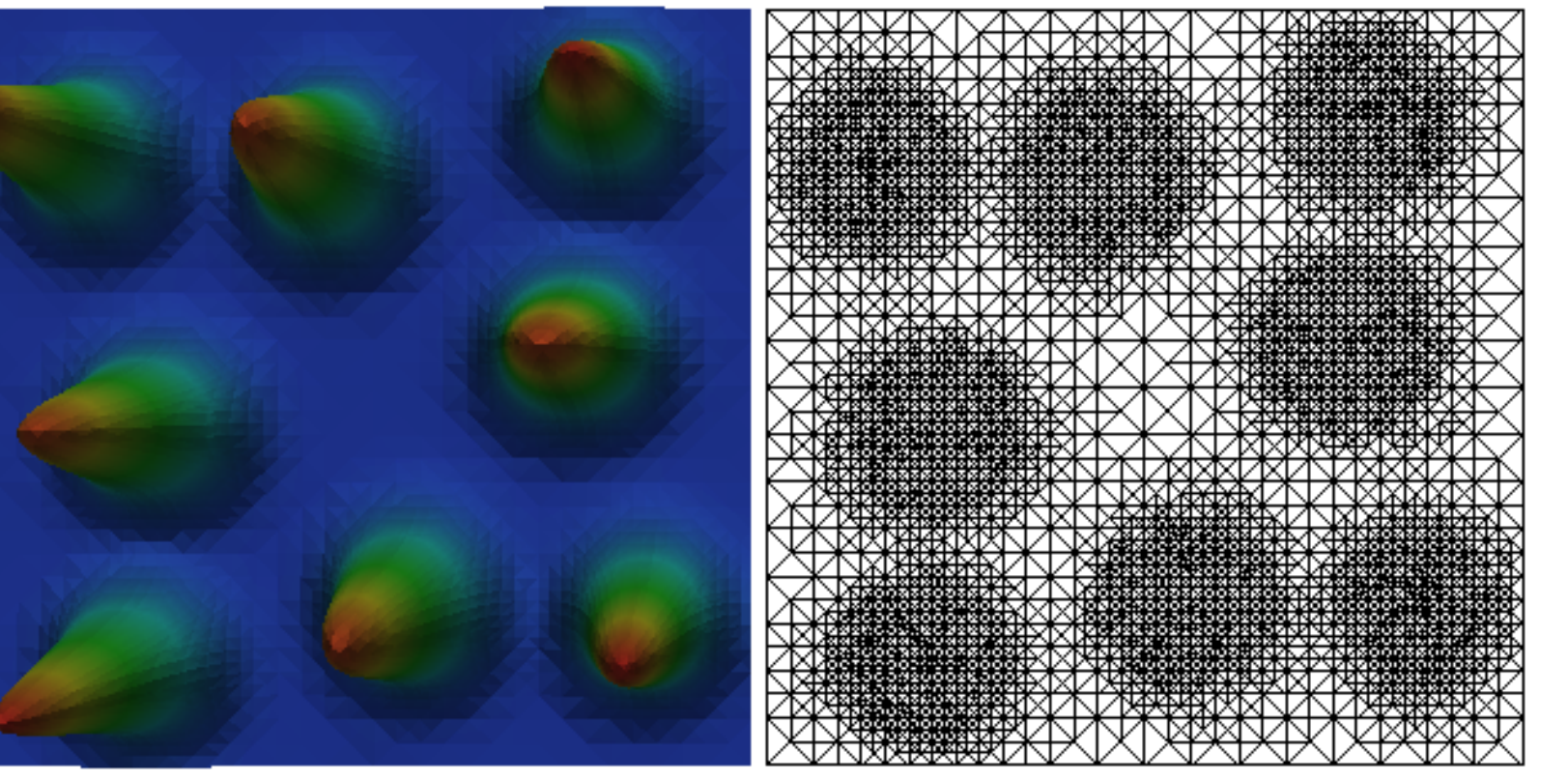}
\label{C:numerics:fig:uniform_t_380}}
\centering
\subfigure[ ][$t=700$]{
\includegraphics[trim = 25mm 185mm 20mm 30mm,  clip, width=0.3\textwidth]{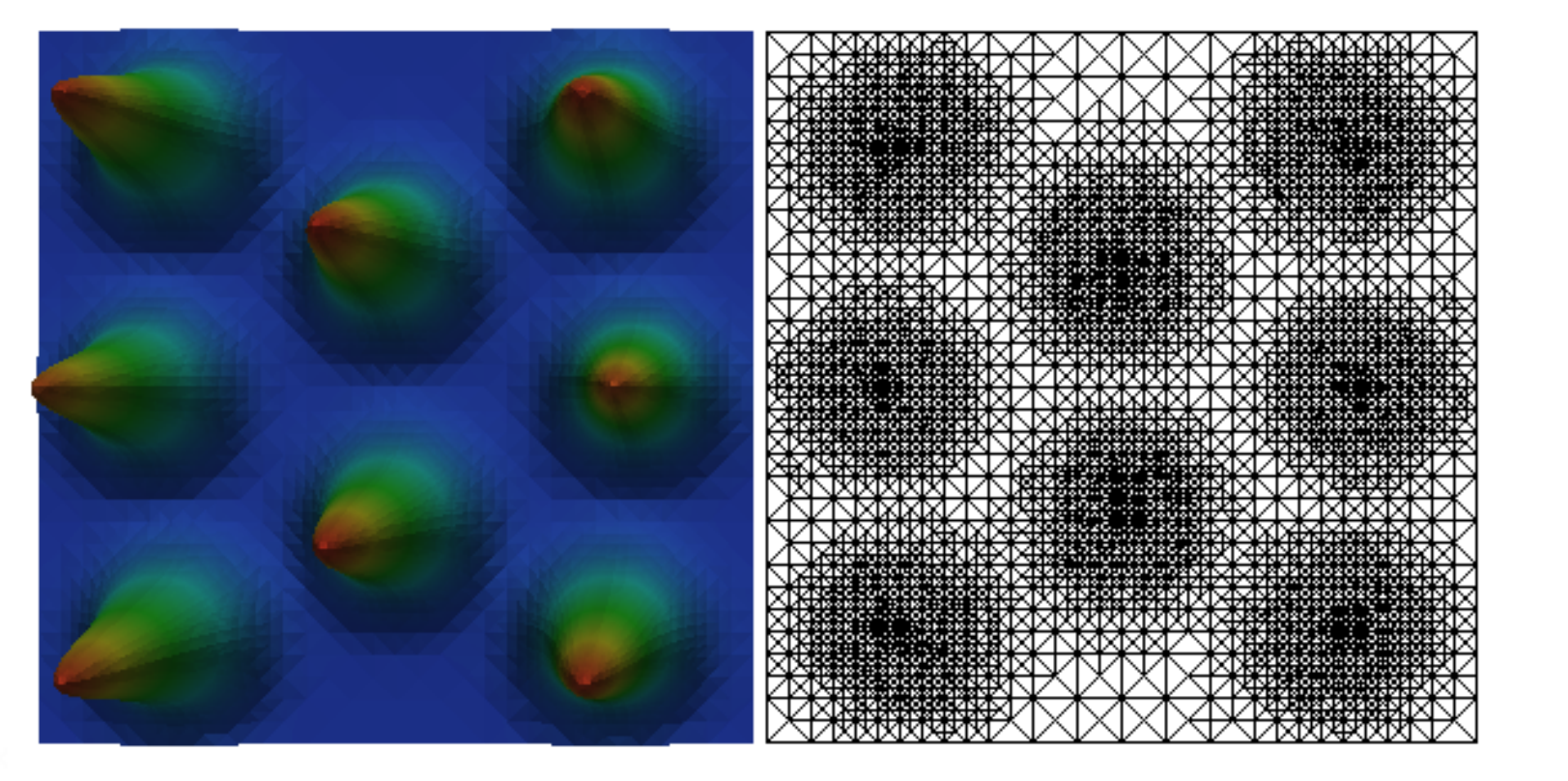}
\label{C:numerics:fig:uniform_t_700}}
\centering
\subfigure[ ][$t=820$]{
\includegraphics[trim = 25mm 185mm 20mm 30mm,  clip, width=0.3\textwidth]{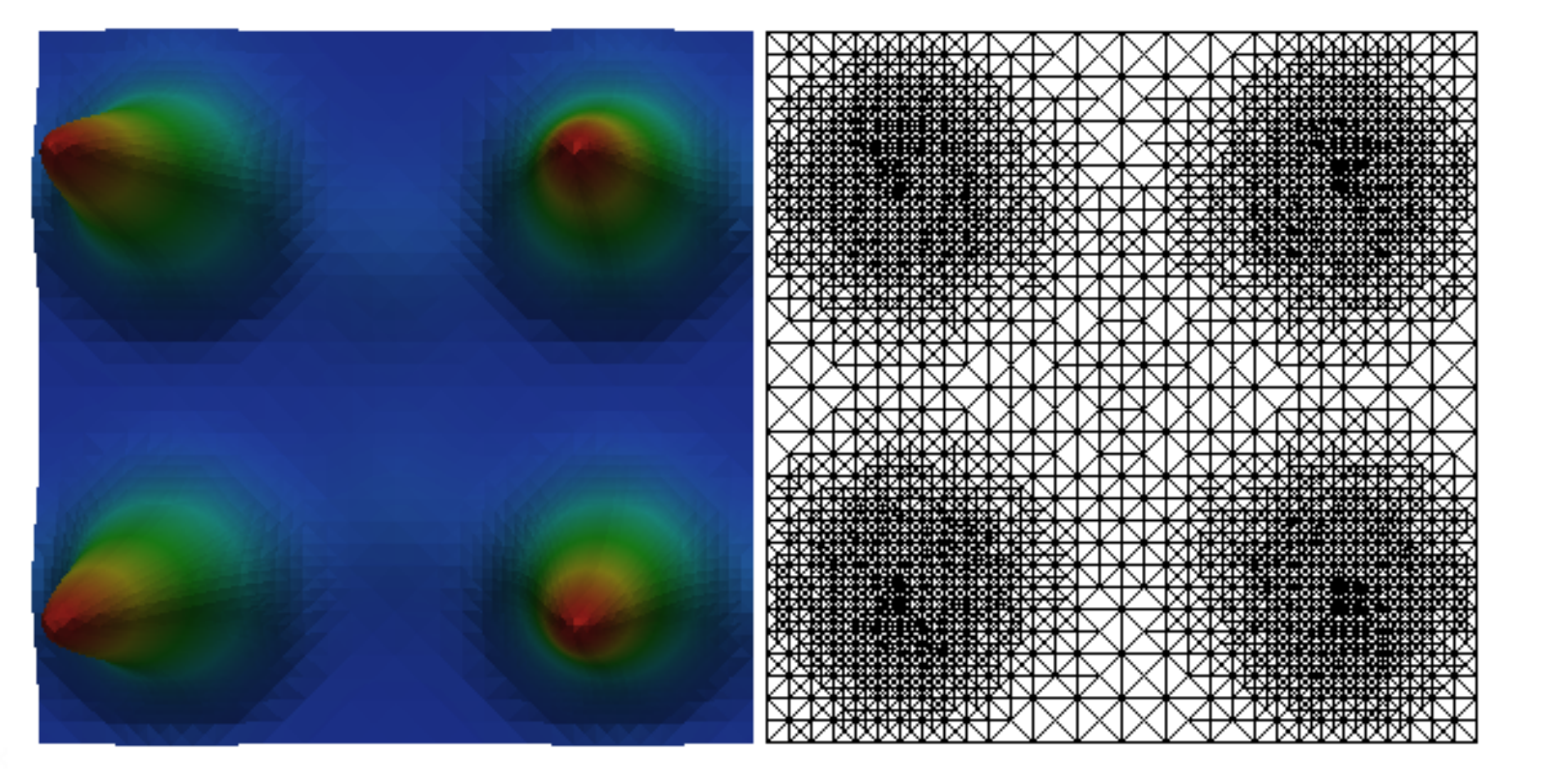}
\label{C:numerics:fig:uniform_t_820}}
\centering
\subfigure[][$t=1000$]{
\includegraphics[trim = 25mm 185mm 20mm 30mm, ,  clip, width=0.3\textwidth]{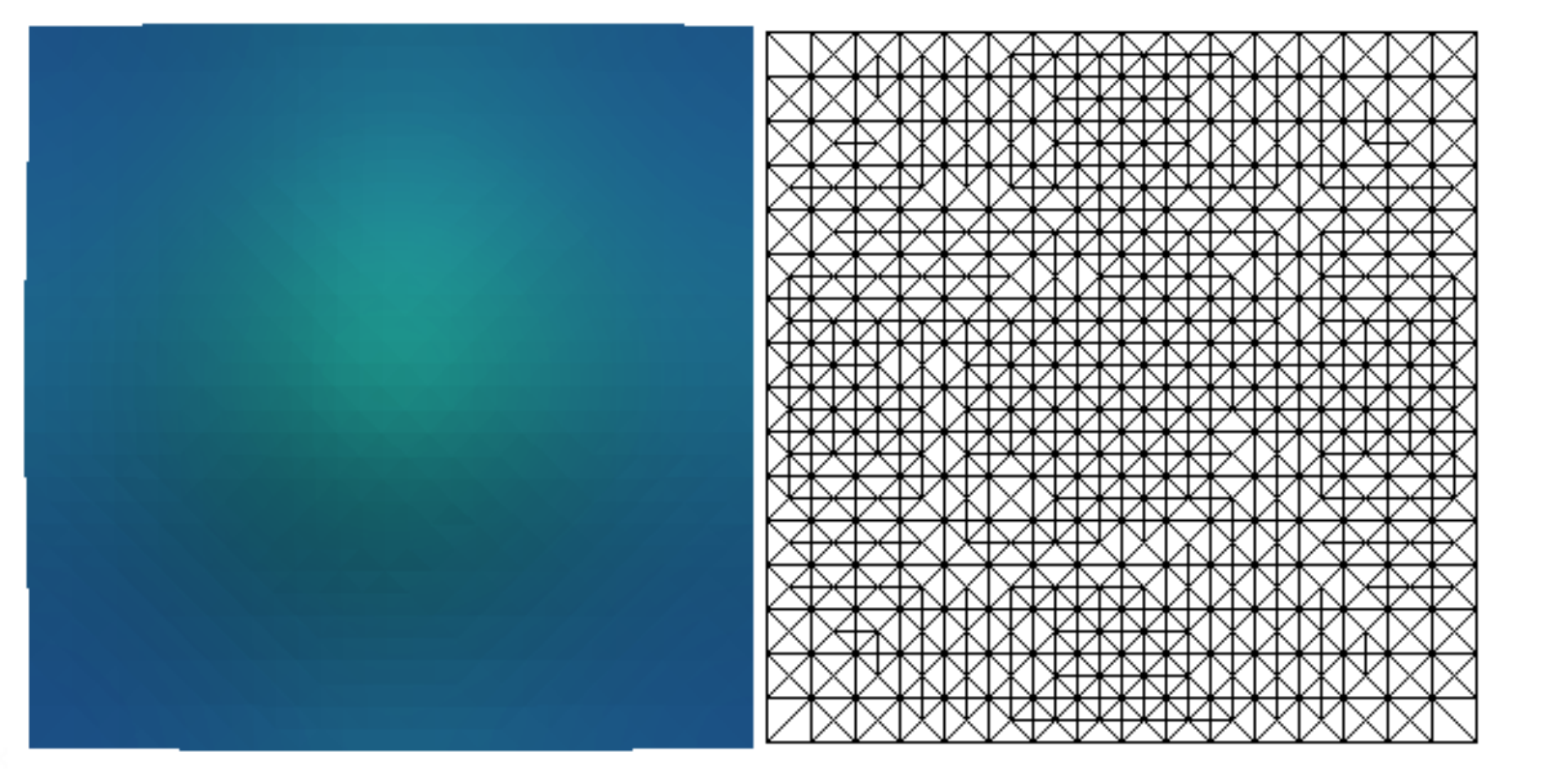}
\label{C:numerics:fig:uniform_t_1000}
}
\caption[]{Snapshots of the discrete activator $(u_1)$ profile for the Schnakenberg kinetics  on the reference domain, under adaptive mesh refinement and domain evolution of the form $(\ref{eqn:SBLI_evolution})$. }  \label{C:numerics:fig:sinusoidal_adaptivity}
\end{figure}

We next present results for the Schnakenberg kinetics, with parameter values $\vec D=(0.01,1)^\transpose$,  $k_1=0.1$,  $k_2=0.9$,  $\gamma=.1,$ where no exact solution is known on a domain with evolution of the form
\begin{equation}\label{eqn:SBLI_evolution}
\A(\vec \xi,t)=\vec \xi(1+9\sin(\pi t/1000)), \quad\vec \xi\in[0,1]^2\text{ and }t\in[0,1000]. 
\end{equation}
We consider an adaptive scheme based on the equidistribution marking strategy with parameters $\theta=0.8$, $tol=10^{-4}$ and a fixed timestep of $10^{-2}$.
Figures \ref{C:numerics:fig:periodic_DOF} and \ref{C:numerics:fig:sinusoidal_adaptivity} show the evolution of the degrees of freedom (DOFs) and the change in discrete solution and snapshots of the activator profiles (on the reference domain) respectively.  The number of DOFs appears positively correlated with the domain size. The mesh is also well refined around the patterns during the evolution, illustrating the benefits of adaptive mesh refinement.

We finish with an application to the case where the evolving domain is an {\it evolving surface} embedded in $\Reals^3$ that is diffeomorphic to a time-independent planar domain. We have derived the model equations and corresponding finite element method elsewhere \cite{amago} and thus only briefly state the details. The model for an RDS posed on an evolving surface is of the form: 
 find $u_i$, functions from $\Gamma_t$ into $\Reals$, such that \fim $u_i$ satisfies
 \begin{equation}\label{C:derivation:eqn:surface_model_problem}
 \begin{split}
\begin{cases}
\pdt u_i(\vec{x},t)+[\vec{a}\cdot\nabla u_i+u_i\nabla_{\Gamma_t}\cdot\vec{a}](\vec{x},t)-
D_i\lap_{\Gamma_t}{u}_i(\vec{x},t)=f_{i}\left(\vec{u}(\vec{x},t)\right), &\vec{x}\in\Gamma_t,t\in(0,T],\\
[\normal\cdot\nabla_{\Gamma_t}{u}_i](\vec{x},t)=0, &\text{ } \vec{x}\in\partial\Gamma_t, t>0,\\
u_i(\vec{x},0)=u_i^{0}(\vec{x}), &\text{ }\vec{x}\in\Gamma_0,
\end{cases}
\end{split}
\end{equation}
 here the Cartesian gradient and Laplacian that appear in (\ref{C:derivation:eqn:model_problem}) are replaced by the surface (tangential) gradient and Laplace-Beltrami operator. 
 We assume the surface  $\Gamma_t$ admits an {\it orthogonal parameterisation} to a planar domain which we denote by 
\begin{equation}
\A:\Oc\subset\Reals^2\times[0,T]\to\Gamma_t\subset\Reals^3.
\end{equation}

%%%%%%%%%%%%%%%%%%%%%%%%%%%%%%min%%%%%%%%%%%%%%%%%%%%%%%%%
%Estimator benchmarking
%%%%%%%%%%%%%%%%%%%%%%%%%%%%%%%%%%%%%%%%%%%%%%%%%%%%%%%%
%%%%%%%%%%%%%%%%%%%%%%%%%%%%%%%%%%%%%%%%%%%%%%%%%%%%%%%
%SPACE ADAPTIVE
%%%%%%%%%%%%%%%%%%%%%%%%%%%%%%%%%%%%%%%%%%%%%%%%%%%%%%%%
%%%%%%%%%%%%%%%%%%%%%%%%%%%%%%%%%%%%%%%%%%%%%%%%%%%%%%%

\begin{figure}[p]
\centering
\subfigure[][$t=50$ (surface)]{
\includegraphics[trim = 0mm 0mm 0mm 0mm,  clip, width=0.225\textwidth]{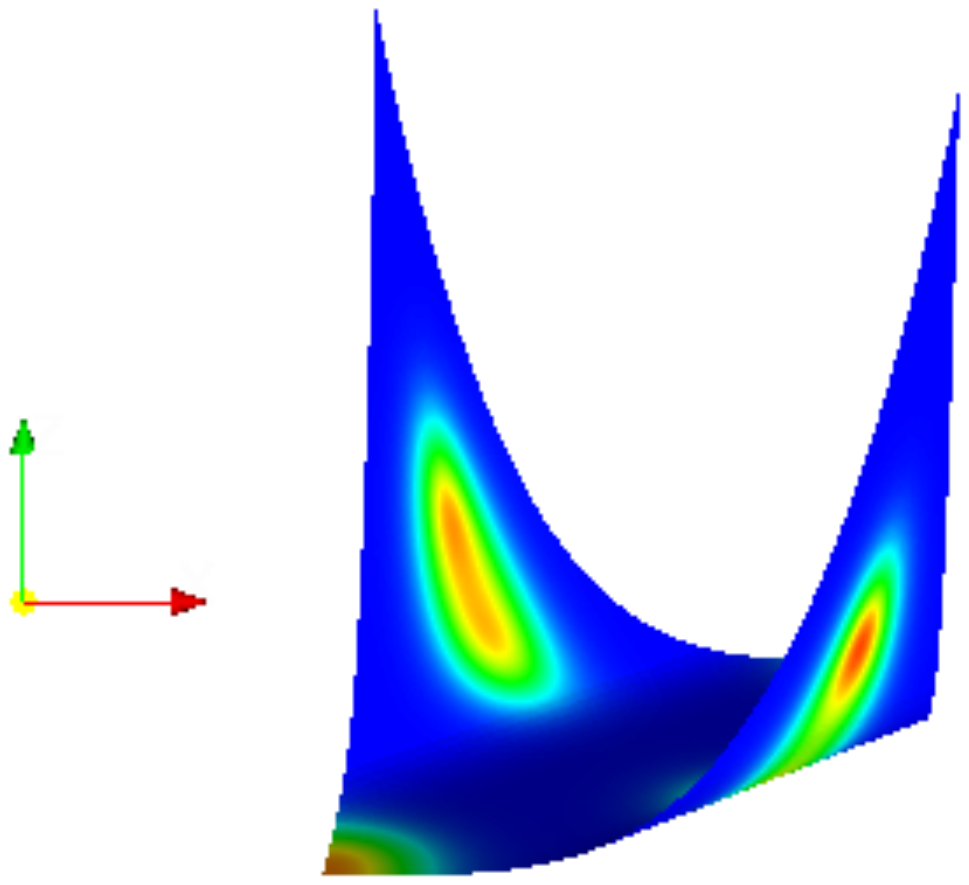}
\label{fig:surface_10}
}
\subfigure[][$t=150$ (surface)]{
%\includegraphcs[trim = 20mm 120mm 20mm 20mm,  clip, width=0.2\textwidth]{./figures/surface/30.pdf}
\includegraphics[trim = 0mm 0mm 0mm 0mm,  clip, width=0.225\textwidth]{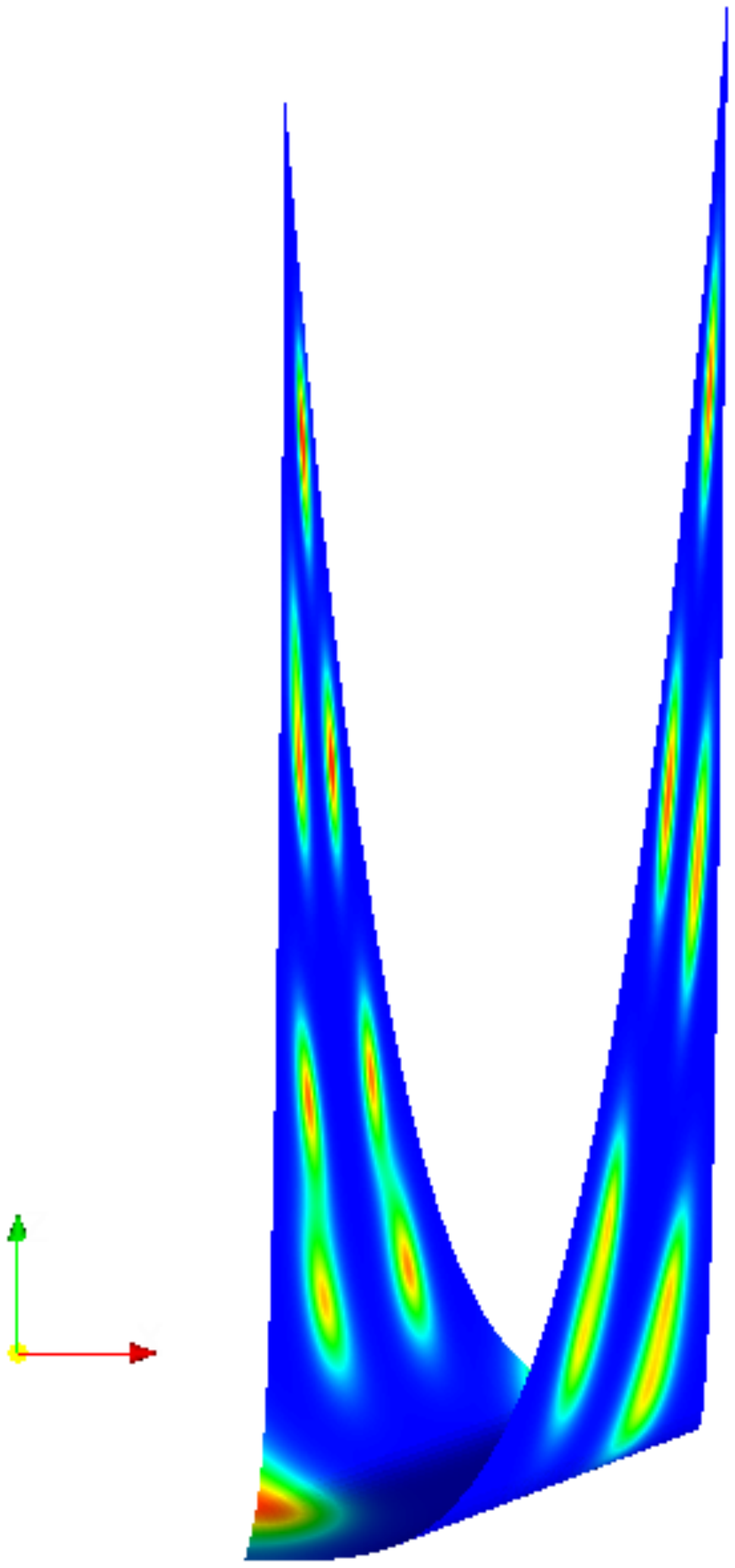}
\label{fig:surface_30}
}
\subfigure[][$t=250$ (surface)]{
\includegraphics[trim = 0mm 0mm 0mm 0mm,  clip, width=0.225\textwidth]{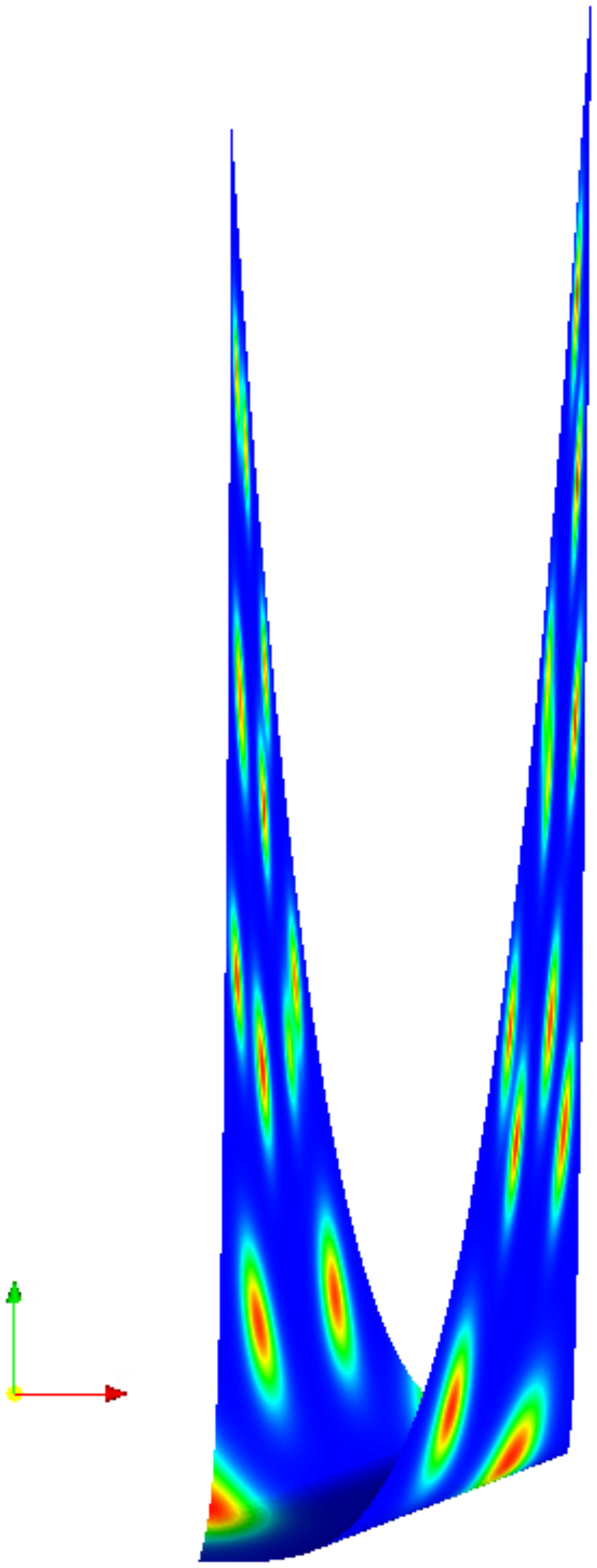}
\label{fig:surface_50}
}
\subfigure[][ $t=50$ (reference and mesh)]{
\includegraphics[trim = 0mm 170mm 0mm 0mm,  clip, width=0.25\textwidth]{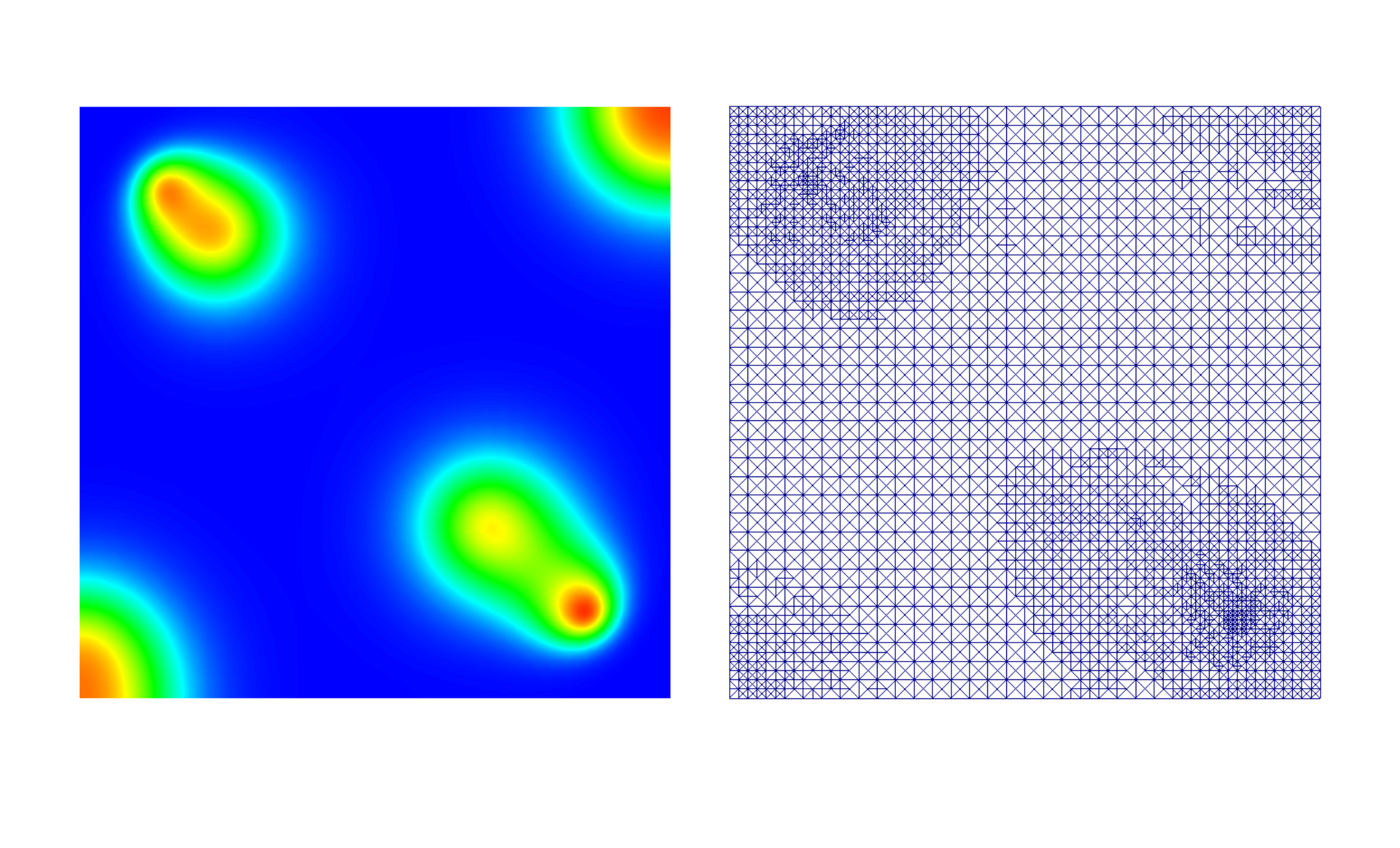}
\label{fig:ref_surface_10}
}
\subfigure[][ $t=150$ (reference and mesh)]{
\includegraphics[trim = 0mm 170mm 0mm 0mm,  clip, width=0.25\textwidth]{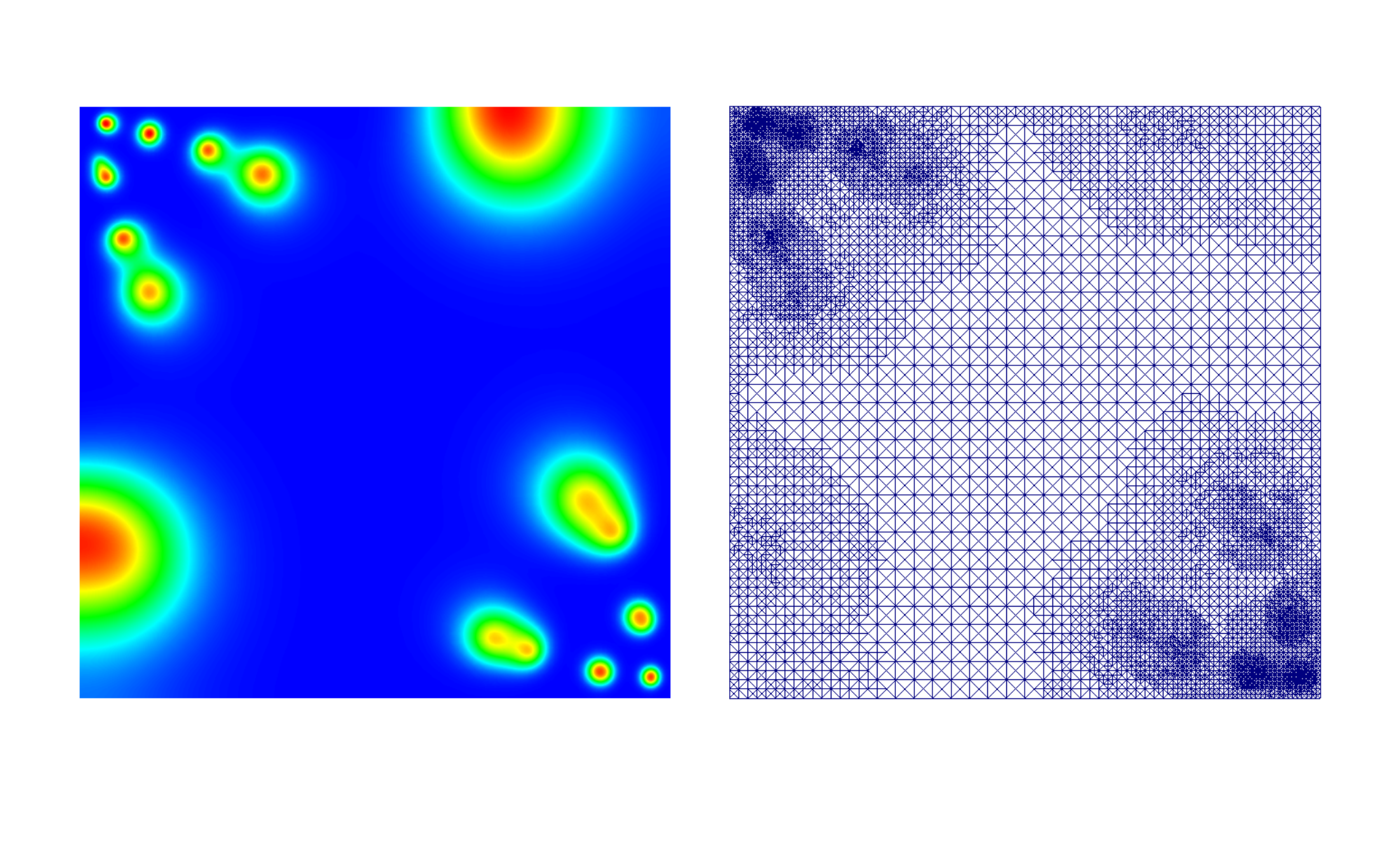}
\label{fig:ref_surface_30}
}
\subfigure[][ $t=250$ (reference and mesh)]{
\includegraphics[trim = 0mm 170mm 0mm 0mm,  clip, width=0.25\textwidth]{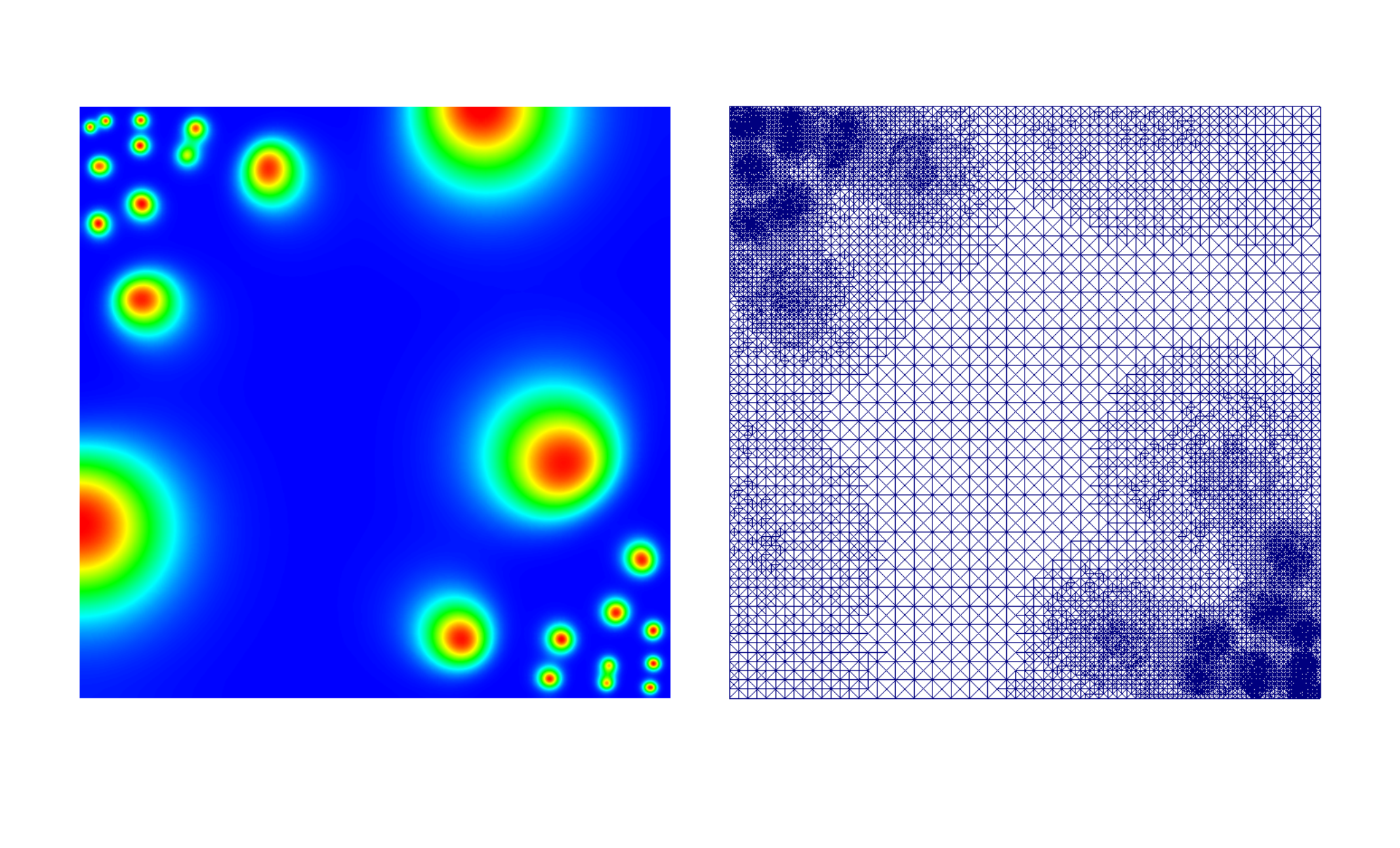}
\label{fig:ref_surface_50}
}
\subfigure[][$t=350$ (surface)]{
\includegraphics[trim = 0mm 0mm 0mm 40mm,  clip, width=0.225\textwidth]{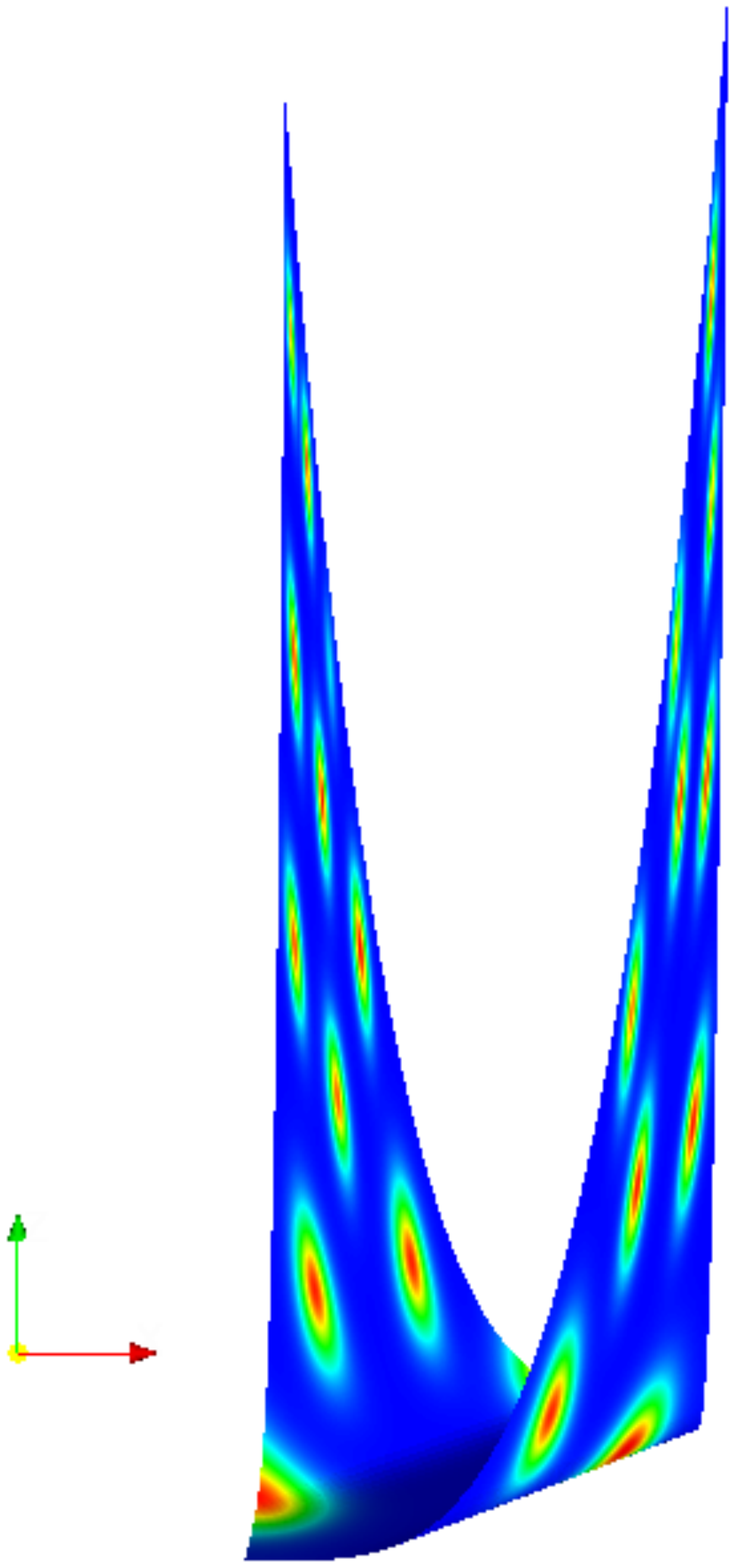}
\label{fig:surface_70}
}
\subfigure[][$t=450$ (surface)]{
\includegraphics[trim = 0mm 0mm 0mm 40mm,  clip, width=0.225\textwidth]{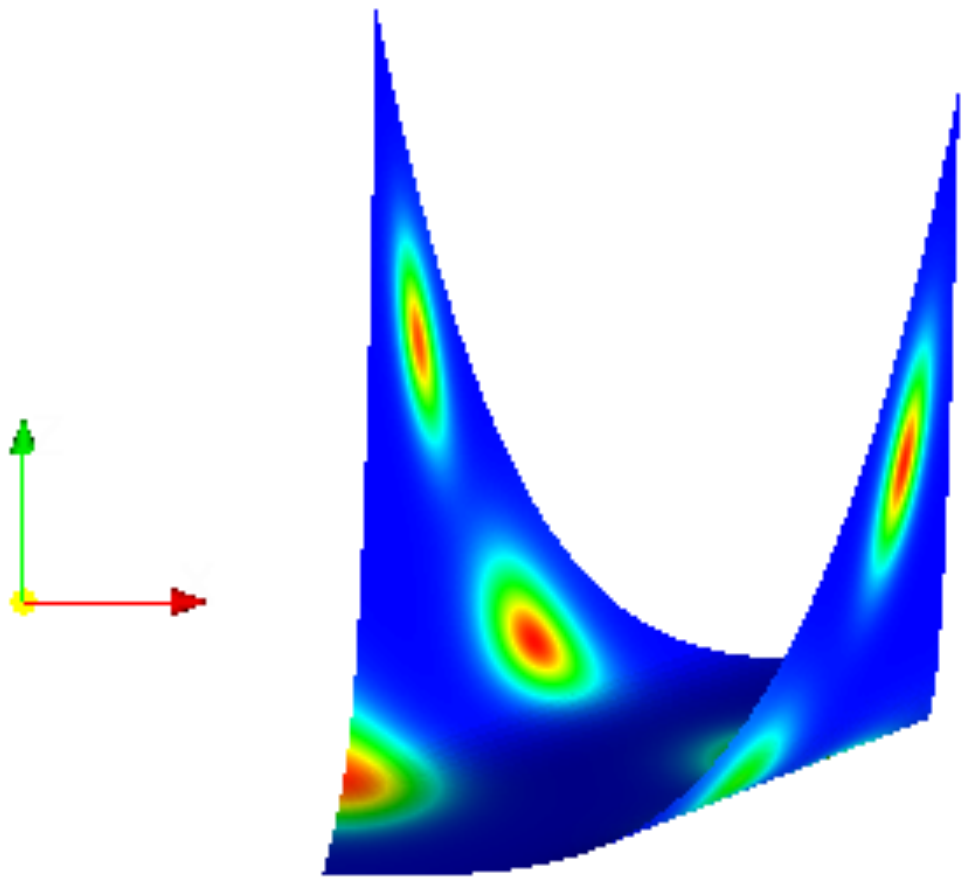}
\label{fig:surface_90}
}
\subfigure[][$t=500$ (surface)]{
\includegraphics[trim = 20mm 120mm 20mm 20mm,  clip, width=0.225\textwidth]{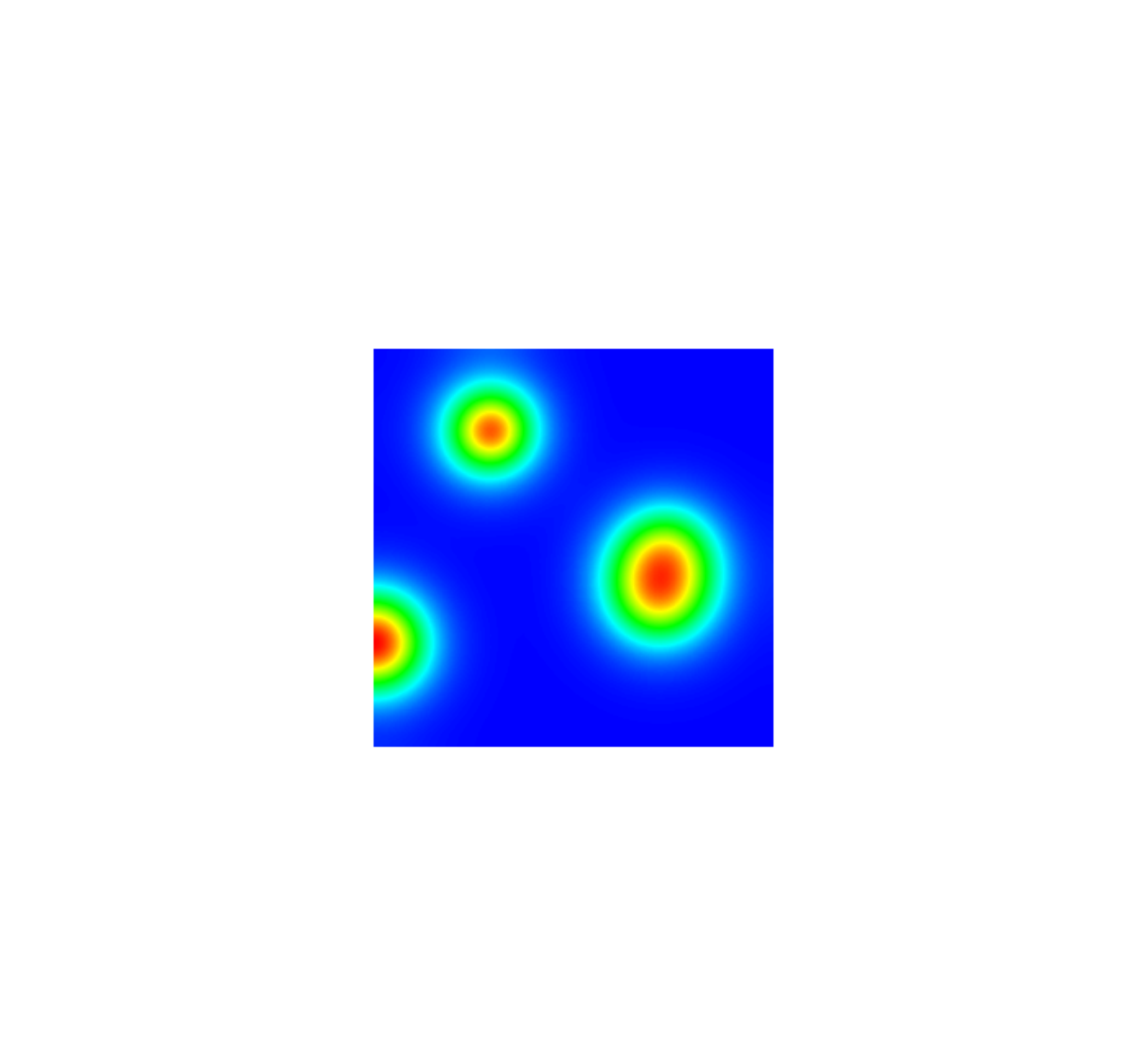}
\label{fig:surface_100}
}
\subfigure[][ $t=350$ (reference and mesh)]{
\includegraphics[trim = 0mm 170mm 0mm 0mm,  clip, width=0.25\textwidth]{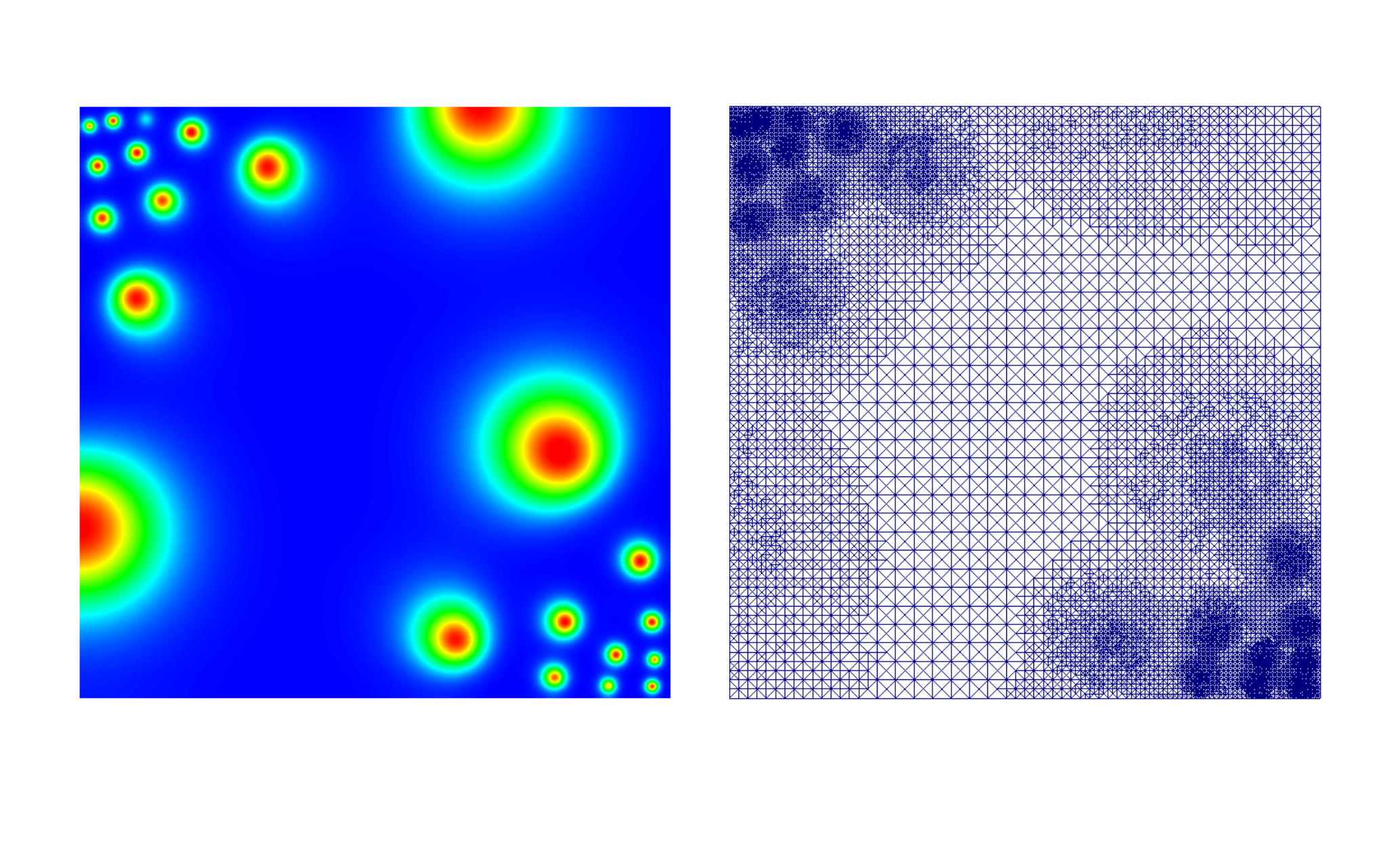}
\label{fig:ref_surface_70}
}
\subfigure[][ $t=450$ (reference and mesh)]{
\includegraphics[trim = 0mm 170mm 0mm 0mm,  clip, width=0.25\textwidth]{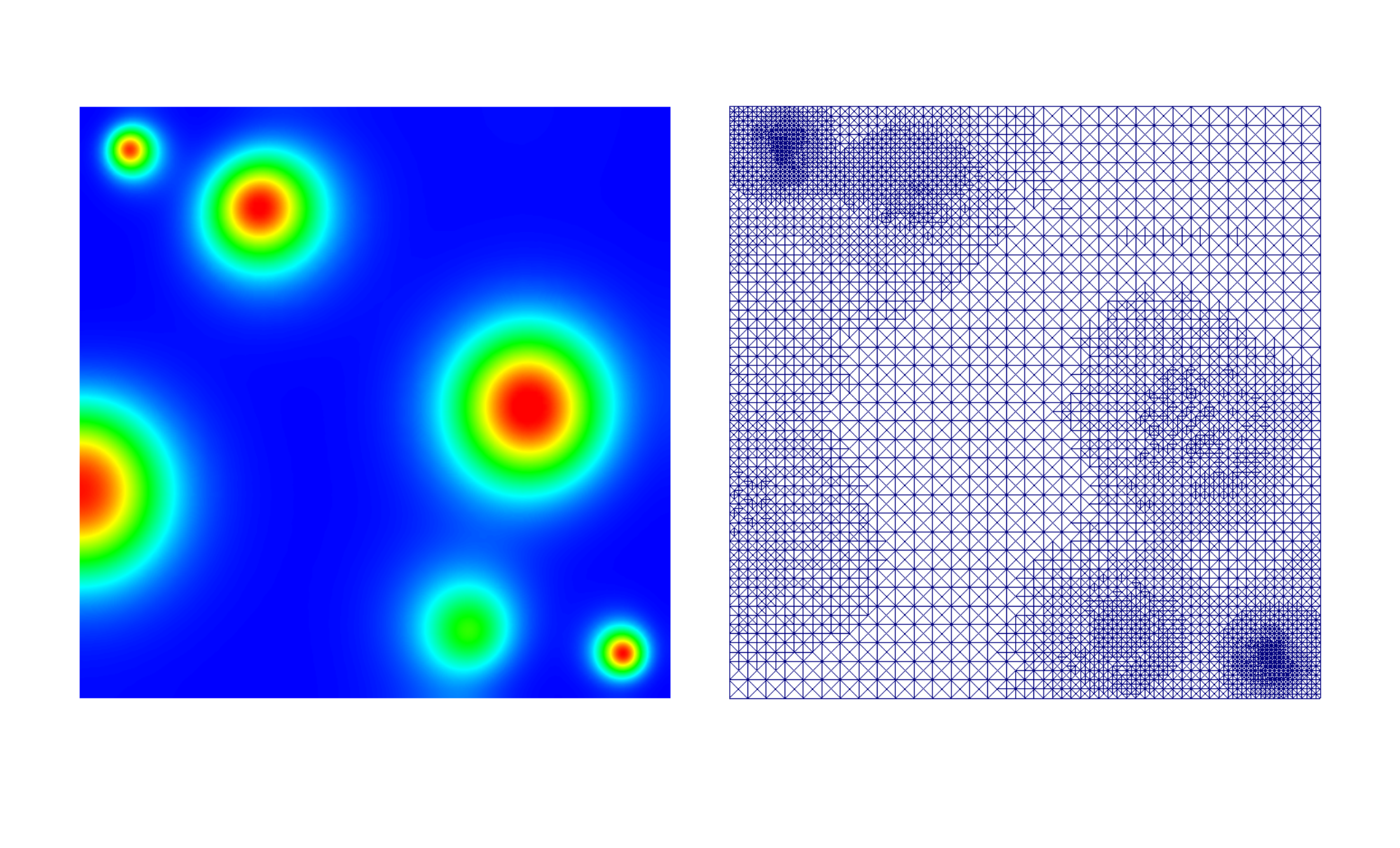}
\label{fig:ref_surface_90}
}
\subfigure[][ $t=500$ (reference and mesh)]{
\includegraphics[trim = 0mm 170mm 0mm 0mm,  clip, width=0.25\textwidth]{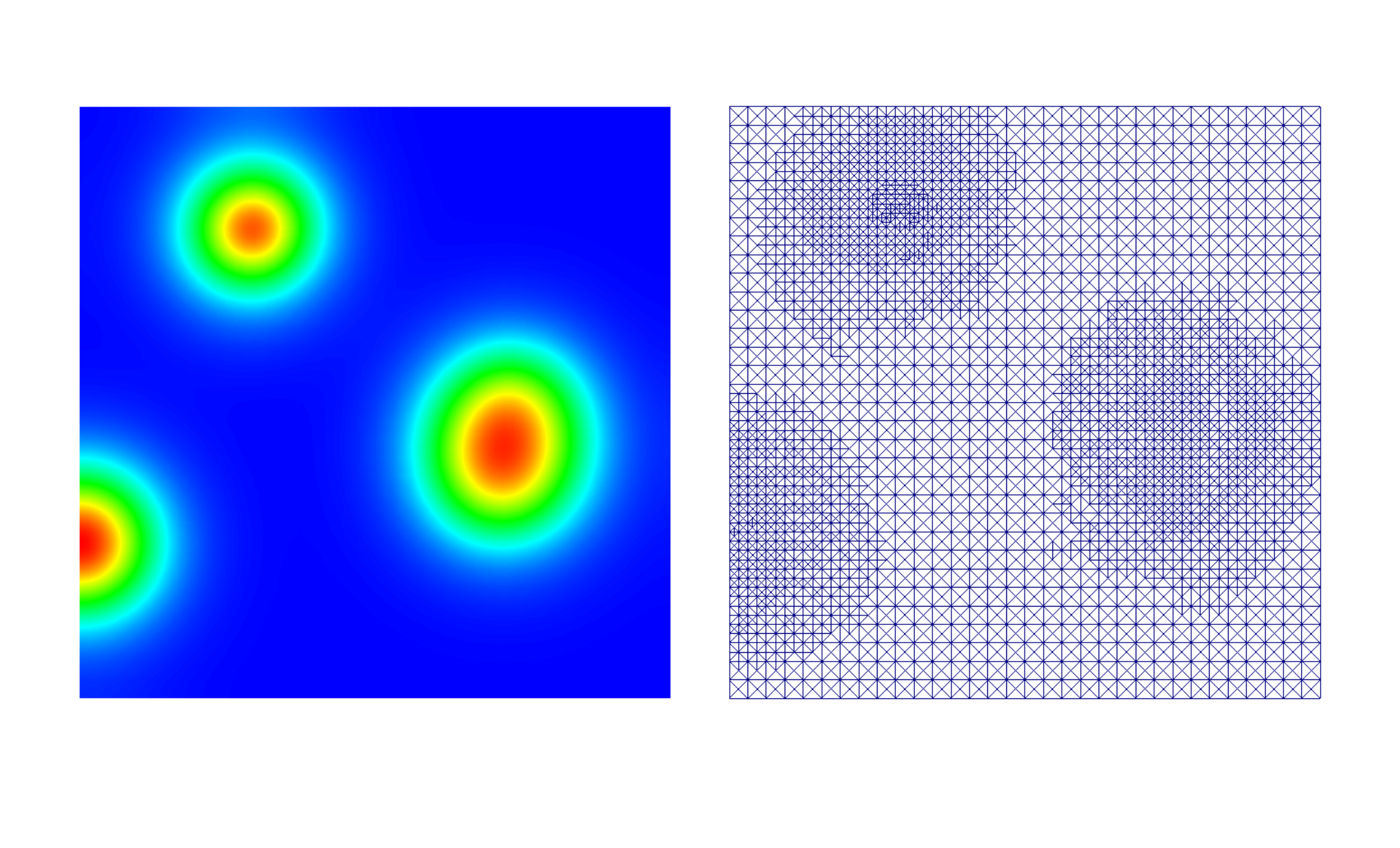}
\label{fig:ref_surface_100}
}
\caption[]{ Snapshots of the discrete activator $(u_1)$ profile for the Schnakenberg kinetics  on the evolving surface  and on the reference domain together with the mesh,  under adaptive mesh refinement and domain evolution of the form (\ref{eqn:surface_evolution}).}  \label{fig:surface_adaptivity}
\end{figure}

Under similar Assumptions to those made in the planar case (see \cite{amago} for details) the corresponding weak formulation on the reference domain $\Oc$ is given by (\ref{C:fem_apriori:eqn:cont_weak_form_reference_ddtju}) where 
the matrix $\vec{K}$  and the determinant of the Jacobian $J$  are given by 
\begin{align}\label{C:applications:eqn:JKdefn}
\vec{K}=\left[\begin{array}{cc}1/{\lv\partial_{1}\A\rv}&0\\0&{1}/{\lv\partial_{2}\A\rv}\end{array}\right]\quad
\text{and} \quad J=\lv\partial_{1}\A\rv\lv\partial_{2}\A\rv.
\end{align}
We consider an example with the Schnakenberg kinetics (\ref{eqn:Schnakenberg_RDS})  , with parameter values $\vec D=(0.01,1)^\transpose$, $k_1=0.1$, $k_2=0.9$, $\gamma=1,$ where no exact solution is known on a domain with evolution of the form
\begin{equation}\label{eqn:surface_evolution}
\A_1( \vec \xi,t)=\xi_1, \ \A_2( \vec \xi,t)=\xi_2, \ \A_3(\vec  \xi,t)=4\sin({\pi t}/{500})(\xi_1-\xi_2)^4\quad\vec \xi\in[0,1]^2,t\in[0,500]. 
\end{equation}
We once again consider an adaptive scheme based on the equidistribution marking strategy with parameters $\theta=0.8$, $tol=10^{-3}$ and a fixed timestep of $10^{-2}$.
Figure \ref{fig:surface_adaptivity} shows snapshots of the activator profiles (on the surface and on the reference domain) and the mesh of the reference domain. As the surface evolves, we observe the emergence of a large number of spots with small radii in the top left and bottom right hand corners of the domain (where curvature is large and growth is fastest)   with annihilation of these spots as the domain contracts. The results clearly illustrate the influence of growth and curvature on pattern formation. The adaptive scheme appears to resolve the solution profiles and the mesh is well refined around the spots on the reference domain, capturing both the small radii spots that develop in the Northwest and Southeast  corners and the large radii spots that develop elsewhere.  

Finally, we remark that we have also considered space-time adaptive schemes based on an heuristic error indicator for the time adaptivity which appear to give dramatic improvements in efficiency \cite{venkataramanthesis}.  
%%%%%%%%%%%%%%%%%%%%%%%%%%%%%%%%%%%%%%%%%%%%%%%%%%%%%%%
%SURFACE
%%%%%%%%%%%%%%%%%%%%%%%%%%%%%%%%%%%%%%%%%%%%%%%%%%%%%%%%
%%%%%%%%%%%%%%%%%%%%%%%%%%%%%%%%%%%%%%%%%%%%%%%%%%%%%%%
%%%%%%%%%%%%%%%%%%%%%%%%%%%%%%%%%%%%%%%%%%%%%%%%%%%%%%%%


\small
\begin{thebibliography}{10}

\bibitem{venkataraman2013adaptive}
C~Venkataraman, O~Lakkis, and A~Madzvamuse.
\newblock Adaptive finite elements for semilinear reaction-diffusion systems on
  growing domains.
\newblock In {\em Numerical Mathematics and Advanced Applications 2011:
  Proceedings of ENUMATH 2011, the 9th European Conference on Numerical
  Mathematics and Advanced Applications, Leicester, September 2011}, page~71.
  Springer, 2013.

\bibitem{ano2000}
A.~Madzvamuse.
\newblock {\em {A Numerical Approach to the Study of Spatial Pattern
  Formation}}.
\newblock PhD thesis, University of Oxford, 2000.

\bibitem{venkataramanthesis}
C.~{Venkataraman}.
\newblock {\em Reaction-diffusion systems on evolving domains with applications
  to the theory of biological pattern formation}.
\newblock PhD thesis, University of Sussex, June 2011.

\bibitem{venk10fem}
O.~{Lakkis}, A.~{Madzvamuse}, and C.~{Venkataraman}.
\newblock {Implicit-explicit timestepping with finite element approximation of
  reaction-diffusion systems on evolving domains}.
\newblock {\em ArXiv e-prints ({\it In Press SINUM})}, 2013.

\bibitem{kruger2003posteriori}
O.~Kruger, M.~Picasso, and JF~Scheid.
\newblock {A posteriori error estimates and adaptive finite elements for a
  nonlinear parabolic problem related to solidification}.
\newblock {\em {Computer Methods in Applied Mechanics and Engineering}},
  {192}({5-6}):{535--558}, 2003.

\bibitem{venkataraman2010global}
Chandrasekhar Venkataraman, Omar Lakkis, and Anotida Madzvamuse.
\newblock Global existence for semilinear reaction--diffusion systems on
  evolving domains.
\newblock {\em Journal of Mathematical Biology}, 64:41--67, 2012.
\newblock 10.1007/s00285-011-0404-x.

\bibitem{medina1996error}
J.~Medina, M.~Picasso, and J.~Rappaz.
\newblock {Error estimates and adaptive finite elements for nonlinear
  diffusion-convection problems}.
\newblock {\em Mathematical Models and Methods in Applied Sciences},
  6(5):689--712, 1996.

\bibitem{madzvamuse2007modified}
A.~Madzvamuse.
\newblock A modified backward euler scheme for advection-reaction-diffusion
  systems.
\newblock {\em Mathematical Modeling of Biological Systems, Volume I}, pages
  183--189, 2007.

\bibitem{schmidt2005design}
A.~Schmidt and K.G. Siebert.
\newblock {\em {Design of adaptive finite element software: The finite element
  toolbox ALBERTA}}.
\newblock Springer Verlag, 2005.

\bibitem{prigo68}
R.~Lefever and I.~Prigogine.
\newblock {Symmetry-breaking instabilities in dissipative systems II}.
\newblock {\em J. chem. Phys}, 48:1695--1700, 1968.

\bibitem{amago}
C.~Venkataraman, T.~Sekimura, E.A. Gaffney, P.K. Maini, and A.~Madzvamuse.
\newblock Modeling parr-mark pattern formation during the early development of
  amago trout.
\newblock {\em Phys. Rev. E}, 84:041923, Oct 2011.

\end{thebibliography}
\end{document}